# The Optimal Momentum Map

Juan-Pablo Ortega
Tudor S. Ratiu

*Dedicated to Jerry Marsden on the occasion of his 60th birthday*

ABSTRACT  The presence of symmetries in a Hamiltonian system usually implies the existence of conservation laws that are represented mathematically in terms of the dynamical preservation of the level sets of a *momentum mapping*. The symplectic or Marsden–Weinstein reduction procedure takes advantage of this and associates to the original system a new Hamiltonian system with fewer degrees of freedom. However, in a large number of situations, this standard approach does not work or is not efficient enough, in the sense that it does not use all the information encoded in the symmetry of the system. In this work, a new momentum map will be defined that is capable of overcoming most of the problems encountered in the traditional approach.

## Contents

# 1   Introduction

Let $(M, \omega, h)$ be a Hamiltonian system and $G$ be a Lie group with Lie algebra $\mathfrak{g}$, acting canonically on $M$; $\omega$ denotes the symplectic two-form on the phase space $M$ and $h : M \to \mathbb{R}$ is the Hamiltonian function. The triplet $(M, \omega, h)$ is called a $G$–**Hamiltonian system** or one says that $(M, \omega, h)$ has **symmetry** $G$, if $h$ is a $G$–invariant function. The $G$–action on $M$ is said to be **globally Hamiltonian** if there exists a $G$–equivariant map $\mathbf{J} : M \to \mathfrak{g}^*$ with respect to the $G$–action on $M$ and the coadjoint action on the dual $\mathfrak{g}^*$ of the Lie algebra $\mathfrak{g}$, such that, for each $\xi \in \mathfrak{g}$ the vector field associated to the infinitesimal generator $\xi_M$ is Hamiltonian with Hamiltonian function $\mathbf{J}^\xi := \langle \mathbf{J}, \xi \rangle$ (the symbol $\langle \cdot, \cdot \rangle$ denotes the natural pairing of $\mathfrak{g}$ with its dual $\mathfrak{g}^*$). The map $\mathbf{J}$ is called the **momentum map** associated to the canonical $G$–action on $M$. The main interest in finding the symmetries of a given system lies in the conservation laws associated to them provided by the following classical result due to E. Noether (see Noether [1918]).

**1.1 Theorem** (Noether). *Let $(M, \omega, h)$ be a $G$–Hamiltonian system. If the $G$–action on $M$ is globally Hamiltonian with associated momentum map $\mathbf{J} : M \to \mathfrak{g}^*$, then $\mathbf{J}$ is a constant of the motion for $h$, that is:*

$$\mathbf{J} \circ F_t = \mathbf{J},$$

*where $F_t$ is the flow of $X_h$, the Hamiltonian vector field associated to $h$.*

In other words, for each $\mu \in \mathfrak{g}_{\mathbf{J}}^* := \mathbf{J}(M)$, the (connected components of the) level set $\mathbf{J}^{-1}(\mu)$ is preserved by the dynamics induced by any $G$–invariant Hamiltonian. Notice that this allows us to look at $\mathfrak{g}_{\mathbf{J}}^*$ as a set of labels that index a family of sets that are invariant under the flows associated to $G$–invariant Hamiltonian functions. The problem with this *classical* approach to the interplay between symmetries and conservation laws resides in the fact that in a number of important situations it cannot be implemented or, even if it can be implemented, it is grossly inefficient in the sense that the sets labeled by $\mathfrak{g}_{\mathbf{J}}^*$ are not the smallest subsets of $M$ preserved by $G$–invariant dynamics. The following situations exemplify these problems:

(i) The simplest situation in which the labeling by $\mathfrak{g}_{\mathbf{J}}^*$ is not optimal is when **the level sets $\mathbf{J}^{-1}(\mu)$ are not connected**. Notice that $G$–invariant dynamics preserves not only $\mathbf{J}^{-1}(\mu)$, that is, the sets labeled by $\mathfrak{g}_{\mathbf{J}}^*$, but also their connected components. Even though in several important situations (for instance, canonical representations of compact connected groups [Lerman, 1995, Theorem 2.1]) the level sets of the momentum map are connected, this is not the case in general.

(ii) **Singularities and the law of conservation of the isotropy**. Let $m \in M$ be such that $G_m \neq \{e\}$. If $\dim G_m > 0$ then $\mathbf{J}(m) = \mu \in \mathfrak{g}^*$ is a singular value of **J** since $\operatorname{rank}(T_m \mathbf{J}) = (\mathfrak{g}_m)^\circ$, where $\mathfrak{g}_m \neq \{0\}$ is the Lie algebra of the isotropy subgroup $G_m$ of $m \in M$ and $(\mathfrak{g}_m)^\circ$ denotes the annihilator of $\mathfrak{g}_m$ in $\mathfrak{g}^*$. Suppose now that the $G$–action on $M$ is proper and denote $H := G_m$. In that situation, the connected components of the set

$$M_H = \{z \in M \mid G_z = H\},$$

are embedded symplectic submanifolds of $(M, \omega)$. Moreover, given that the flow of any $G$–invariant Hamiltonian is $G$–equivariant, the connected components of $M_H$ are preserved by any of these flows (law of conservation of the isotropy), that is, in the presence of points with non trivial symmetry, the smallest sets left invariant by $G$–invariant dynamics are not the connected components of $\mathbf{J}^{-1}(\mu)$ but rather those of $\mathbf{J}^{-1}(\mu) \cap M_H$.

**(iii) Symmetries given by finite groups**. This is an important particular case of **(ii)**. Many relevant systems possess symmetries that are expressed through the canonical action of a finite group $G$ on $M$. Since in that case the Lie algebra $\mathfrak{g} = \{0\}$, the associated momentum map **J** is the constant map equal to zero. Therefore, in this scenario, the Noether Theorem is empty of content.

**(iv) Symmetries without a momentum map**. As the statement of Noether's Theorem implies, one needs to insure a priori the existence of a momentum map which then gives the conserved quantities associated to the given canonical symmetry. However, even if the system possesses a canonical symmetry, the existence of a momentum map is by no means guaranteed. For example, let $M = S^1 \times S^1 = \mathbb{T}^2$ with the symplectic form $\omega = d\theta_1 \wedge d\theta_2$. Let $G = S^1$ acting on $M$ by $e^{i\phi} \cdot (e^{i\theta_1}, e^{i\theta_2}) := (e^{i(\phi+\theta_1)}, e^{i\theta_2})$. A simple calculation shows (see Weinstein [1976]) that this action is canonical but that it does not admit an associated globally defined momentum map.

The conservation laws of a Hamiltonian system allows one to apply **symplectic** or **Marsden–Weinstein reduction** (Marsden and Weinstein [1974]; Meyer [1973]) in order to obtain a new system with less degrees of freedom. However, in some of the above mentioned problematic situations, performing reduction becomes either impossible or a task subject to arbitrary choices.

The main goal of this paper is the construction of a momentum map, which we will call **optimal momentum map**, that does not suffer from the inconveniences pointed out in the classical approach. This mapping is always defined in the presence of a canonical symmetry, it uses the symmetries of the system in order to provide optimal conservation laws, and can be used, following the traditional Marsden–Weinstein scheme, to symplectically reduce the system in virtually every possible situation.

We emphasize that this paper is just an introduction to the optimal momentum map. Many of its properties are still being investigated. For instance, the convexity properties of its image will be presented in a future publication. Some natural questions about this object have been recently answered; for instance how to carry out orbit reduction and reduction by stages in the optimal framework is already well understood (Marsden, Misiolek, Ortega, Perlmutter, and Ratiu [2001]; Ortega [2001a]). Incidentally, this shows how to perform standard Hamiltonian singular reduction by stages without using the so called *stages hypothesis*.

The paper is organized as follows. Section 2 introduces some notations and technical results that will be needed later on. The definition of the optimal momentum map is given in Section 3. This then leads to a Noether Theorem valid also in the above mentioned problematic situations. Section 4 shows how the optimal momentum map can be used to perform Marsden–Weinstein reduction. In addition, it is shown that these new reduced spaces coincide with the usual ones, both in the regular (Marsden and Weinstein [1974]) and singular (Sjamaar and Lerman [1991]; Bates and Lerman [1997]; Ortega [1998]; Ortega and Ratiu [2002]) situations.

## 2  Preliminaries and technical results

### 2.1  Generalized distributions

We begin by recalling some results on generalized distributions due to Stefan [1974a,b] and Sussman [1973]. The reader is also encouraged to check with the excellent review in the Appendix 3 of the book by Libermann and Marle [1987].

**2.1 Definition.** *Let $M$ be a differentiable manifold. A **generalized distribution** $D$ on $M$ is a subset of the tangent bundle $TM$ such that, for any point $m \in M$, the fiber $D_m = D \cap T_m M$ is a vector subspace of $T_m M$. The dimension of $D_m$ is called the **rank** of the distribution $D$ at $m$. A **differentiable section** of $D$ is a differentiable vector field $X$ defined on an open subset $U$ of $M$, such that for any point $z \in U$, $X(z) \in D_z$. An immersed connected submanifold $N$ of $M$ is said to be an **integral manifold** of the distribution $D$ if, for every $z \in N$, $T_z i(T_z N) \subset D_z$, where $i : N \to M$ is the canonical injection. The integral submanifold $N$ is said to be of **maximal dimension** at a point $z \in N$ if $T_z i(T_z N) = D_z$.*

(i)  *The generalized distribution $D$ is **differentiable** if, for every point $z \in M$, and for every vector $v \in D_z$, there exists a differentiable section $X$ of $D$, defined on an open neighborhood $U$ of $z$, such that $X(z) = v$.*

(ii)  *The generalized distribution $D$ is **completely integrable** if, for every*

point $z \in M$, there exists an integral manifold of $D$ everywhere of maximal dimension which contains $z$.

(iii) The generalized distribution $D$ is **involutive** if it is invariant under the (local) flows associated to differentiable sections of $D$.

**Remark.** Our definition of involutive distribution is more general than the traditional one which states that $D$ is involutive if $[X, Y]$ takes values in $D$ whenever $X$ and $Y$ are vector fields with values in $D$. The two concepts of involutivity are equivalent only when the dimension of $D_m$ is the same for any $m \in M$. ♦

**2.2 Theorem** (Generalized Frobenius Theorem). *A differentiable distribution $D$ on a manifold $M$ is completely integrable iff it is involutive.*

**Proof.** See Stefan [1974a,b]; Sussman [1973]. ∎

In our discussion we will be interested in the specific case in which the generalized distribution is given by an everywhere defined family of vector fields, that is, there is a family of smooth vector fields $\mathcal{D}$ whose elements are vector fields $X$ defined on a open subset $\mathrm{Dom}(X) \subset M$ such that, for any $z \in M$ the generalized distribution $D$ is given by

$$D_z = \mathrm{span}\{X(z) \in T_z M | X \in \mathcal{D} \text{ and } z \in \mathrm{Dom}(X)\}.$$

Note that in such a case the distribution $D$ is differentiable by construction. We will say that $D$ is the generalized distribution **spanned** by $\mathcal{D}$. One of the reasons for our interest in this special case resides in the fact that when these distributions are completely integrable, a very useful characterization of their integral manifolds can be given. In order to describe it we introduce some terminology following Libermann and Marle [1987]. Let $X$ be a vector field defined on an open subset $\mathrm{Dom}(X)$ of $M$ and $F_t$ be its flow. For any fixed $t \in \mathbb{R}$ the domain $\mathrm{Dom}(F_t)$ of $F_t$ is an open subset of $\mathrm{Dom}(X)$ such that $F_t : \mathrm{Dom}(F_t) \to \mathrm{Dom}(F_{-t})$ is a diffeomorphism. If $Y$ is a second vector field defined on the open set $\mathrm{Dom}(Y)$ with flow $G_t$ we can consider, for two fixed values $t_1, t_2 \in \mathbb{R}$, the composition of the two diffeomorphisms $F_{t_1} \circ G_{t_2}$ as defined on the open set $\mathrm{Dom}(G_{t_2}) \cap (G_{t_2})^{-1}(\mathrm{Dom}(F_{t_1}))$ (which may be empty).

The previous prescription allows us to inductively define the composition of an arbitrary number of locally defined flows. We will obviously be interested in the flows associated to the vector fields in $\mathcal{D}$ that define the distribution $D$. The following sentences describe some important conventions that we will use all over the paper. Let $k \in \mathbb{N}^*$, be a positive natural number, $\mathcal{X}$ be an ordered family $\mathcal{X} = (X_1, \ldots, X_k)$ of $k$ elements of $\mathcal{D}$, and $T$ be a $k$–tuple $T = (t_1, \ldots, t_k) \in \mathbb{R}^k$ such that $F_t^i$ denotes the (locally defined) flow of $X_i$, $i \in \{1, \ldots, k\}$, $t_i$. We will denote by $\mathcal{F}_T$ the locally defined diffeomorphism $\mathcal{F}_T = F_{t_1}^1 \circ F_{t_2}^2 \circ \cdots \circ F_{t_k}^k$ constructed using the above given prescription. Any diffeomorphism from an open subset of $M$ onto another open subset of $M$ that is constructed in the same fashion as $\mathcal{F}_T$ is said to be **generated** by the family $\mathcal{D}$. It can be proven that the composition of diffeomorphisms generated by $\mathcal{D}$ and the inverses of diffeomorphisms generated by $\mathcal{D}$ are themselves diffeomorphisms generated by $\mathcal{D}$ [Libermann and Marle, 1987, Proposition 3.3, Appendix 3]. In other words, the family of diffeomorphisms generated by $\mathcal{D}$ forms a **pseudogroup of transformations** (see page 74 of Paterson [1999]) that will be denoted by $G_\mathcal{D}$. Two points $x$ and $y$ in $M$ are said to be $G_\mathcal{D}$–**equivalent**, if there exists a diffeomorphism $\mathcal{F}_T \in G_\mathcal{D}$ such that $\mathcal{F}_T(x) = y$. The relation *being $G_\mathcal{D}$–equivalent* is an equivalence relation whose equivalence classes are called the $G_\mathcal{D}$–**orbits**.

**2.3 Theorem.** *Let $D$ be a differentiable generalized distribution on the smooth manifold $M$ spanned by a family of vector fields $\mathcal{D} \subset \mathfrak{X}(M)$. The following properties are equivalent:*

(i) *The distribution $D$ is invariant under the pseudogroup of transformations generated by $\mathcal{D}$, that is, for each $\mathcal{F}_T \in G_\mathcal{D}$ generated by $\mathcal{D}$ and for each $z \in M$ in the domain of $\mathcal{F}_T$,*

$$T_z \mathcal{F}_T(D_z) = D_{\mathcal{F}_T(z)}.$$

(ii) *The distribution $D$ is completely integrable and its integral manifolds are the $G_\mathcal{D}$–orbits.*

**Proof.** See Stefan [1974a,b]; Sussman [1973]. For a compact presentation combine Theorems 3.9 and 3.10 in the Appendix 3 of Libermann and Marle [1987]. ∎

## 2.2 Poisson manifolds and Poisson tensors

A **Poisson manifold** is a pair $(M, \{\cdot, \cdot\})$, where $M$ is a differentiable manifold and $\{\cdot, \cdot\}$ is a bilinear operation on $C^\infty(M)$ such that $(C^\infty(M), \{\cdot, \cdot\})$ is a Lie algebra, and $\{\cdot, \cdot\}$ is a derivation (that is, the Leibniz identity holds) in each argument; $(C^\infty(M), \{\cdot, \cdot\})$ is called a **Poisson algebra**. The derivation property of the Poisson bracket defines for each $h \in C^\infty(M)$ the **Hamiltonian vector field** $X_h$ by $\mathbf{d}f(X_h) = \{f, h\}$ for any $f \in C^\infty(M)$. Obviously, any Hamiltonian system on a symplectic manifold is a Poisson dynamical system relative to the Poisson bracket induced by the symplectic structure. The converse relation is given by the **Symplectic Stratification Theorem**, (see Kirillov [1976], Weinstein [1983], Libermann and Marle [1987], or Marsden and Ratiu [1999]) which states that any Poisson manifold $(M, \{\cdot, \cdot\})$ is partitioned into **symplectic leaves**. Each leaf is, by definition, the set of points that can be linked to a given one by a finite number of smooth curves, each of which is a piece of an integral curve of a locally defined Hamiltonian vector field. The symplectic leaves are connected immersed symplectic manifolds in $M$ (relative to the inclusion map), whose Poisson bracket coincides with that of $M$. The tangent space at $z \in M$ to a leaf consists of all vectors that are equal to the value of some Hamiltonian vector field at $z$. The symplectic leaves are invariant under the flow of any Hamiltonian vector field on $M$.

The derivation property of the Poisson bracket implies that for any two functions $f, g \in C^\infty(M)$, the value of the bracket $\{f, g\}(z)$ at an arbitrary point $z \in M$ (and therefore $X_f(z)$ as well), depends on $f$ only through $\mathbf{d}f(z)$ (see [Abraham, Marsden, and Ratiu, 1988, Theorem 4.2.16] for a justification of this argument) which allows us to define a contravariant antisymmetric two–tensor $B \in \Lambda^2(T^*M)$ by

$$B(z)(\alpha_z, \beta_z) = \{f, g\}(z),$$

where $\mathbf{d}f(z) = \alpha_z$ and $\mathbf{d}g(z) = \beta_z \in T_z^*M$. This tensor is called the **Poisson tensor** of $M$. We will denote by $B^\sharp : T^*M \to TM$ the vector bundle map associated to $B$, that is

$$B(z)(\alpha_z, \beta_z) = \langle \alpha_z, B^\sharp(\beta_z) \rangle.$$

The Poisson tensor permits another formulation of the results regarding symplectic leaves in terms of the **characteristic distribution**. Given a Poisson manifold $(M, \{\cdot, \cdot\})$ with associated Poisson tensor $B$, the characteristic distribution $D$ on $M$ is defined by $D := B^\sharp(T^*M)$. It can be proven [Libermann and Marle, 1987, Theorem 12.1, Chapter III] that the characteristic distribution $D$ is differentiable, completely integrable, and that its integral manifolds are the symplectic leaves of $(M, \{\cdot, \cdot\})$.

**2.4 Proposition.** *Let $(M, \omega)$ be a symplectic manifold and $B \in \Lambda^2(T^*M)$ be the associated Poisson tensor. Then for any $m \in M$ and any vector subspace $V \subset T_mM$,*

$$V^\omega = B^\sharp(m)(V^\circ),$$

*where $V^\omega := \{v \in T_mM \mid \omega(m)(u, v) = 0, \text{ for all } u \in V\}$ is the $\omega$– **orthogonal complement** of $V$ in $T_mM$.*

**Proof.** Let $\alpha_m \in V^\circ$ and $f \in C^\infty(M)$ be such that $\mathbf{d}f(m) = \alpha_m$. Then, for any $v \in V$:

$$\omega(m)(B^\sharp(m)(\alpha_m), v) = \omega(m)(X_f(m), v) = \mathbf{d}f(m) \cdot v = \langle \alpha_m, v \rangle = 0,$$

which proves that $B^\sharp(m)(V^\circ) \subset V^\omega$. Given that in the symplectic case $B^\sharp(m)$ is an isomorphism for all $m \in M$, a dimension count concludes the proof. ∎

**2.5 Definition.**   *A smooth mapping $\varphi : M_1 \to M_2$, between the two Poisson manifolds $(M_1, \{\cdot,\cdot\}_1)$ and $(M_2, \{\cdot,\cdot\}_2)$ is called **canonical** or **Poisson** if for all $g, h \in C^\infty(M_2)$ we have*

$$\varphi^*\{g, h\}_2 = \{\varphi^*g, \varphi^*h\}_1.$$

For future reference we state a result whose proof can be found, for instance, in [Marsden and Ratiu, 1999, Proposition 10.3.2].

**2.6 Proposition.**   *Let $\varphi : M_1 \to M_2$ be a smooth map between two Poisson manifolds $(M_1, \{\cdot,\cdot\}_1)$ and $(M_2, \{\cdot,\cdot\}_2)$. Then $\varphi$ is a Poisson map if and only if $T\varphi \circ X_{h\circ\varphi} = X_h \circ \varphi$ for any $h \in C^\infty(M_2)$. In particular, if $h \in C^\infty(M_2)$, $F_t^2$ is the flow of $X_h$, and $F_t^1$ is the flow of $X_{h\circ\varphi}$, then $F_t^2 \circ \varphi = \varphi \circ F_t^1$.*

### 2.3   Proper actions, tubes, and slices

The following definitions and results are standard in Lie theory (see Bredon [1972]; Palais [1961]). In what follows we will deal mostly with proper actions. Recall that the action $\Phi : G \times M \to M$ is called **proper** if for any two convergent sequences $\{m_n\}$ and $\{g_n \cdot m_n\}$ in $M$, there exists a convergent subsequence $\{g_{n_k}\}$ in $G$. Proper actions have compact isotropy subgroups and Hausdorff orbit spaces.

**2.7 Definition.**   *Let $G$ be a Lie group and $H \subset G$ be a subgroup. Suppose that $H$ acts on the left on the manifold $A$. The **twist action** of $H$ on the product $G \times A$ is defined by*

$$h \cdot (g, a) = (gh, h^{-1} \cdot a).$$

*Note that this action is free and proper by the freeness and properness of the action on the $G$–factor. The **twisted product** $G \times_H A$ is defined as the orbit space corresponding to the twist action.*

**2.8 Proposition.**   *The twisted product $G \times_H A$ is a $G$–space relative to the left action defined by $g' \cdot [g, a] := [g'g, a]$. If the action of $H$ on $A$ is proper, the $G$–action on $G \times_H A$ just defined is proper.*

**2.9 Definition.**   *Let $M$ be a manifold and $G$ be a Lie group acting on $M$. Let $m \in M$ and denote $H := G_m$. A **tube** about the orbit $G \cdot m$ is a $G$–equivariant diffeomorphism*

$$\varphi : G \times_H A \longrightarrow U,$$

*where $U$ is a $G$–invariant neighborhood of $G \cdot m$ in $M$ and $A$ is some open ball centered at the origin in a representation space of $H$.*

Note that if the $G$–action on $M$ is proper then the $G$–action on any tube $G \times_H A$ is also proper since the isotropy subgroup $H$ is compact and, consequently, its action on $A$ is proper. Proposition 2.8 guarantees the claim.

**2.10 Definition.** *Let $M$ be a manifold and $G$ be a Lie group acting properly on $M$. Let $m \in M$ and denote $H := G_m$. Let $S$ be a submanifold of $M$, such that $m \in S$ and $H \cdot S = S$. We say that $S$ is a **slice** at $m$ if the map*

$$\begin{array}{rcl} G \times_H S & \longrightarrow & U \\ [g, s] & \longmapsto & g \cdot s \end{array}$$

*is a tube about $G \cdot m$, for some $G$–invariant open neighborhood $U$ of $G \cdot m$.*

**2.11 Theorem.** *Let $M$ be a manifold and $G$ be a Lie group acting properly on $M$. Let $m \in M$, denote $H := G_m$, and let $S$ be a submanifold of $M$ such that $m \in S$. Then the following statements are equivalent:*

(i) *There is a tube $\varphi : G \times_H A \longrightarrow U$ about $G \cdot m$ such that $\varphi[e, A] = S$.*

(ii) *$S$ is a slice at $m$.*

(iii) *$G \cdot S$ is an open neighborhood of $G \cdot m$ and there is an equivariant retraction*

$$r : G \cdot S \longrightarrow G \cdot m$$

*such that $r^{-1}(m) = S$.*

The ball $A$ appearing in (i) can always be chosen to be a $H$–invariant neighborhood of $0$ in the vector space $T_m M / T_m(G \cdot m)$, where $H$ acts linearly and orthogonally by $h \cdot [v] = [h \cdot v]$.

**2.12 Theorem** (Slice Theorem)**.** *Let $M$ be a manifold and $G$ be a Lie group acting properly on $M$. Then there is a slice for the $G$–action at any $m \in M$.*

As a first consequence of the Slice Theorem we have the following result that will be used in the sequel.

**2.13 Proposition.** *Let $G$ be a Lie group acting properly on the manifold $M$, $U$ an open $G$–invariant subset of $M$, and $f \in C^\infty(U)^G$. Then for any $z \in U$ there exist a $G$–invariant open neighborhood $V \subset U$ of $z$ and a $G$–invariant smooth function $F \in C^\infty(M)^G$ such that $f|_V = F|_V$.*

**Proof.** Let $U_1 \subset U$ be an open $G$–invariant neighborhood of $z$ that by the Slice Theorem can be modeled by the tube $U_1 \simeq G \times_{G_z} A_r$, where $A_r$ is the open ball of radius $r$ in the vector space $T_z M / T_z(G \cdot z)$ on which $G_z$ acts orthogonally. Define $g : A_r \to \mathbb{R}$ as the smooth $G_z$–invariant function given by $g(v) := f([e, v])$. Since $G_z$ is compact, there exists a $G_z$–invariant bump function $\phi : A_r \to [0, 1]$ such that

$$\phi|_{D_{r/2}} = 1 \quad \text{and} \quad \phi|_{D_r \setminus D_{3r/4}} = 0.$$

Define $f' \in C^\infty(U_1)^G$ by $f'([h, v]) := \phi(v)g(v)$, for any $h \in G$ and $v \in A_r$. Since $f'$ and all its derivatives vanish on the boundary of $U_1$, its extension off $U_1$ by the identically zero function yields $F \in C^\infty(M)^G$. Take $V \simeq G \times_{G_z} A_{r/2}$. It is clear that $F|_V = f'|_V = f|_V$. ∎

The following result, proved for the first time in Ortega [1998], will be of great importance in the sequel.

**2.14 Proposition.** *Let $G$ be a Lie group acting properly on the smooth manifold $M$. Let $m \in M$ be a point with isotropy subgroup $H := G_m$. Then*

$$((T_m(G \cdot m))^\circ)^H = \mathrm{span}\{\mathbf{d}f(m) \mid f \in C^\infty(M)^G\}.$$

**Proof.** We first show that if $f \in C^\infty(M)^G$, then $\mathbf{d}f(m) \in ((T_m(G\cdot m))^\circ)^H$. It is clear that for any $\xi \in \mathfrak{g}$,

$$\langle \mathbf{d}f(m), \xi_M(m) \rangle = \frac{d}{dt}\bigg|_{t=0} f(\exp t\xi \cdot m) = \frac{d}{dt}\bigg|_{t=0} f(m) = 0.$$

Hence, $\mathbf{d}f(m) \in T_m(G \cdot m)^\circ$. Now, $\mathbf{d}f(m)$ is also $H$–fixed since for any $h \in H$ and any $v = \frac{d}{dt}\big|_{t=0} m(t) \in T_m M$ with $m(0) = m$,

$$\langle h \cdot \mathbf{d}f(m), v \rangle = \langle \mathbf{d}f(m), h^{-1} \cdot v \rangle$$
$$= \frac{d}{dt}\bigg|_{t=0} f(h^{-1} \cdot m(t))$$
$$= \frac{d}{dt}\bigg|_{t=0} f(m(t)) = \langle \mathbf{d}f(m), v \rangle.$$

Since the vector $v$ is arbitrary, $h \cdot \mathbf{d}f(m) = \mathbf{d}f(m)$, as required.

Since we are going to work locally, in order to prove the converse inclusion, we do it in the tubular model provided by the Slice Theorem. Thus, the manifold $M$ will be replaced by $G \times_H V$, where $V = T_m M / T_m(G \cdot m)$ and the point $m \in M$ is represented by $[e, 0] \in G \times_H V$. It is easy to verify that

$$T_{[e, 0]}(G \cdot [e, 0]) = \{T_{(e, 0)}\pi(\xi, 0) \in T_{[e, 0]}(G \times_H V) \mid \xi \in \mathfrak{g}\} \cong \mathfrak{g}/\mathfrak{h} \times \{0\},$$

where $\pi : G \times V \to G \times_H V$ is the canonical projection. Clearly,

$$(T_{[e, 0]}(G \cdot [e, 0]))^\circ \cong \{0\} \times V^* \cong V^*. \tag{2.1}$$

At this point we introduce the following

**2.15 Lemma.** [1] *Let $H$ be a compact Lie group acting linearly on the vector space $V$, as well as on its dual $V^*$ via the associated contragredient representation. Then, the restriction to $(V^*)^H$ of the dual map associated to the inclusion $i_{V^H} : V^H \hookrightarrow V$ is a $H$-equivariant isomorphism from $(V^*)^H$ to $(V^H)^*$.*

**Proof.** Note that for any $\beta \in V^*$, $i_{V^H}^*(\beta) = \beta|_{V^H}$. Take an $H$-invariant inner product $\langle\langle \cdot, \cdot \rangle\rangle$ on $V$, always available by the compactness of $H$. Let $W$ be the $H$-invariant orthogonal complement to $V^H$ with respect to $\langle\langle \cdot, \cdot \rangle\rangle$.

Any element $\alpha \in (V^H)^*$ can be extended to $\beta \in V^*$ by setting $\beta|_W = 0$. Moreover, note that $i_{V^H}^*(\beta) = \beta|_{V^H} = \alpha$ and also, for any $v \in V^H$, $w \in W$, and $h \in H$, we have

$$\langle h \cdot \beta, v + w \rangle = \langle \beta, h^{-1} \cdot (v+w) \rangle = \langle \beta, v + h^{-1} \cdot w \rangle = \langle \beta, v + w \rangle,$$

since both $w$ and $h^{-1} \cdot w$ are in $W$. This implies that $\beta \in (V^*)^H$, and hence $i_{V^H}^*|_{(V^*)^H}$ is surjective.

For injectivity, suppose $\beta|_{V^H} = 0$, for some $\beta \in (V^*)^H$. Let $v \in V$ be such that $\langle \beta, w \rangle = \langle\langle v, w \rangle\rangle$, for all $w \in W$. Then, for any $h \in H$ and any $w \in V$ we have that

$$\langle\langle h \cdot v, w \rangle\rangle = \langle\langle v, h^{-1} \cdot w \rangle\rangle = \langle \beta, h^{-1} \cdot w \rangle = \langle h \cdot \beta, w \rangle = \langle \beta, w \rangle = \langle\langle v, w \rangle\rangle,$$

which, by the non degeneracy of the inner product implies that $h \cdot v = v$ and hence $v \in V^H$. But $\beta|_{V^H} = 0$, which implies that $v = 0$, which in turn implies $\beta = 0$. Hence $i_{V^H}|_{(V^*)^H}$ is an isomorphism.

The $H$-equivariance of $i_{V^H}^*|_{(V^*)^H}$ follows trivially from the following chain of equalities that are satisfied for any $h \in H$, $\beta \in (V^*)^H$, and $v \in V^H$

$$\langle h \cdot i_{V^H}^*(\beta), v \rangle = \langle i_{V^H}^*(\beta), h^{-1} \cdot v \rangle = \langle \beta, h^{-1} \cdot v \rangle = \langle \beta, v \rangle = \langle i_{V^H}^*(\beta), v \rangle.$$

▼

Now, using Lemma 2.15 in (2.1) we get

$$((T_{[e,\,0]}(G \cdot [e,\,0]))^\circ)^H \simeq (V^*)^H \cong (V^H)^*.$$

In the tubular model, the $G$–invariant functions $f \in C^\infty(G \times_H V)^G$ are characterized by the condition $f \circ \pi \in C^\infty(V)^H$. The claim then follows if we show that

$$(V^*)^H = \{\mathbf{d}g(0) \in V^* \mid g \in C^\infty(V)^H\}.$$

---

[1] We thank Tanya Schmah for her quick proof of this lemma.

Let $g \in C^\infty(V)^H$ and $h \in H$ be arbitrary. Then, for any $v = \frac{d}{dt}\big|_{t=0} c(t) \in V$ with $c(0) = 0$, we have

$$\begin{aligned}
\langle h \cdot \mathbf{d}g(0), v \rangle &= \langle \mathbf{d}g(0), h^{-1} \cdot v \rangle \\
&= \frac{d}{dt}\bigg|_{t=0} g(h^{-1} \cdot c(t)) \\
&= \frac{d}{dt}\bigg|_{t=0} g(c(t)) = \langle \mathbf{d}g(0), v \rangle.
\end{aligned}$$

Since $v \in V$ is arbitrary, it follows that $h \cdot \mathbf{d}g(0) = \mathbf{d}g(0)$.

To prove the converse, we begin by decomposing $V$ into its irreducible $H$–components:

$$V = W_1 \oplus \ldots \oplus W_k \oplus U_1 \oplus \ldots \oplus U_r,$$

where $\dim W_1 = \ldots = \dim W_k = 1$, and $\dim U_i > 1$ for $i \in \{1, \ldots, r\}$. Thus,

$$V^H = W_1 \oplus \ldots \oplus W_k. \tag{2.2}$$

Let $\{w_1, \ldots, w_k\}$, be a basis of $V^H$ adapted to the splitting (2.2). Define $\pi_1, \ldots, \pi_k \in V^*$ by

$$\begin{aligned}
\pi_i(w_j) &= \delta_{ij} \quad i, j \in \{1, \ldots, k\} \\
\pi_i|_{U_p} &= 0 \quad i \in \{1, \ldots, k\},\ p \in \{1, \ldots, r\}.
\end{aligned}$$

By construction, the functionals $\pi_1, \ldots, \pi_k \in V^*$ are linear invariants of the $H$–action on $V$. Moreover, they are the only ones. Indeed, since $\pi_1, \ldots, \pi_k$ is a basis of $(V^H)^*$, there are no additional independent linear invariants on $V^H$. If $\alpha : U_1 \oplus \ldots \oplus U_r \to \mathbb{R}$ is another nontrivial linear invariant, there is some $p \in \{1, \ldots, r\}$ such that $\alpha|_{U_p}$ is not the zero functional. Therefore, $\ker(\alpha|_{U_p}) \neq 0$ is a nontrivial $H$–invariant subspace of $U_p$. Since this is impossible by the irreducibility of $U_p$, it follows that such an $\alpha$ cannot exist.

We have thus shown that $\pi_1, \ldots, \pi_k \in V^*$, or, in general, that any basis of $(V^H)^*$ spans the set of all independent linear invariants of the $H$–action on $V$. By the Hilbert Theorem, the ring of $H$–invariant polynomials on $V$ is finitely generated. We complete the set $\{\pi_1, \ldots, \pi_k\}$ to a generating system $\{\pi_1, \ldots, \pi_k, \pi_{k+1}, \ldots, \pi_q\}$ of the this ring. The Schwarz Theorem (Schwarz [1974]; Mather [1977]) guarantees that every $H$–invariant function $f \in C^\infty(V)^H$ can be locally written as

$$f = g(\pi_1, \ldots, \pi_q),$$

with $g \in C^\infty(\mathbb{R}^q)$. Let now $\alpha \in (V^*)^H \cong (V^H)^*$ be arbitrary. The form $\alpha \in (V^H)^*$ can be expanded as

$$\alpha = \alpha_1 \pi_1 + \ldots + \alpha_k \pi_k.$$

with $\alpha_1, \ldots, \alpha_k \in \mathbb{R}$. Let $g \in C^\infty(\mathbb{R}^q)$ be such that

$$\frac{\partial g(0)}{\partial \pi_i} = \alpha_i, \quad i \in \{1, \ldots, k\}.$$

With this choice, the function $f := g(\pi_1, \ldots, \pi_q)$ belongs to $C^\infty(V)^H$ and satisfies

$$\mathbf{d}f(0) = \frac{\partial g(0)}{\partial \pi_1}\pi_1 + \ldots + \frac{\partial g(0)}{\partial \pi_k}\pi_k = \alpha_1\pi_1 + \ldots + \alpha_k\pi_k = \alpha,$$

where we used that $\mathbf{d}\pi_j(0) = 0$ for $j \in \{k+1, \ldots, q\}$ because the invariants $\pi_j$ in this range of the indices are at least quadratic. Since $\alpha$ is arbitrary, the result follows. ∎

**Remark.** The properness condition in the statement of the previous proposition is essential (and is not tied to the existence of slices) since there are examples of non proper actions where this result does not hold. Indeed, consider the irrational flow on the torus. Since the orbits of this action fill densely the torus, the only invariant functions in this particular case are the constant functions. Hence the right hand side of the equality in Proposition 2.14 is trivial. However, if the torus in question is bigger than one dimensional, the vector space $(T_m(G \cdot m))^\circ$ is non trivial. ♦

We now recall a few facts that we will need later on about the interplay between group actions and symplectic and Poisson structures.

**2.16 Definition.** *Let $(M, \{\cdot, \cdot\})$ be a Poisson manifold (resp. $(M, \omega)$ a symplectic manifold), let $G$ be a Lie group, and let $\Phi : G \times M \to M$ be a smooth left action of $G$ on $M$. We say that the action $\Phi$ is **canonical** if $\Phi$ acts by canonical transformations; that is, for any $f, h \in C^\infty(M)$ and any $g \in G$*

$$\Phi_g^*\{f, h\} = \{\Phi_g^*f, \Phi_g^*h\} \qquad (\text{resp. } \Phi_g^*\omega = \omega).$$

**2.17 Proposition.** *Let $(M, \{\cdot, \cdot\})$ be a Poisson manifold and denote by $B \in \Lambda^2(T^*M)$ be the associated Poisson tensor. Let $G$ be a Lie group acting canonically on $M$. Then, for any $m \in M$ such that $G_m =: H$ and any vector subspace $V \subset T_m^*M$:*

(i) *$B^\sharp(m) : T_m^*M \to T_mM$ is $H$–equivariant.*

(ii) *If the Poisson bracket $\{\cdot, \cdot\}$ is induced by a symplectic manifold $\omega$ then*

$$B^\sharp(m)(V^H) = (B^\sharp(m)(V))^H,$$

*where the $H$–superscript denotes the set of $H$–fixed points in the corresponding spaces.*

**Proof.** Part (i) is a trivial consequence of the canonical character of the action. Part (ii) follows from Part (i) and the non degeneracy of $B^\sharp(m)$ in the symplectic case. ∎

## 3 A new momentum map and an optimal Noether Theorem

In this section we introduce the main ideas of the paper.

### 3.1 The optimal momentum map

**3.1 Definition.** *Let $(M, \{\cdot, \cdot\})$ be a Poisson manifold, $G$ be a Lie group acting canonically on $M$, $U$ be a $G$–invariant open subset of $M$, and $C^\infty(M)^G$ (respectively $C^\infty(U)^G$) be the set of $G$–invariant functions on $M$ (respectively, $G$–invariant functions on $U$). Let $\mathcal{E}$ be the set of Hamiltonian vector fields associated to all the elements of $C^\infty(U)^G$, for all the open $G$–invariant subsets $U$ of $M$, that is,*

$$\mathcal{E} = \left\{ X_f \mid f \in C^\infty(U)^G, \text{ with } U \subset M \text{ open and } G\text{–invariant} \right\},$$

*and $E$ be the smooth generalized distribution on $M$ spanned by $\mathcal{E}$. We will call $E$ the $G$–**characteristic distribution**.*

**Remark.** If the $G$–action on $M$ is proper, the definition of the distribution $E$ admits some simplification. Indeed, by definition, for any $m \in M$, there is a $r \in \mathbb{N}$ such that

$$E(m) = \text{span}\{X_{f_1}(m), \ldots, X_{f_r}(m)\},$$

where $f_i \in C^\infty(U_i)^G$ for any $i \in \{1, \ldots, r\}$. By Proposition 2.13, for each $i \in \{1, \ldots, r\}$ there exists an open $G$–invariant subset $V_i \subset U_i$ containing the point $m$ and a function $F_i \in C^\infty(M)^G$ such that $f_i|_{V_i} = F_i|_{V_i}$. Consequently,

$$E(m) = \text{span}\{X_{f_1}(m), \ldots, X_{f_r}(m)\} = \text{span}\{X_{F_1}(m), \ldots, X_{F_r}(m)\}.$$

This proves that the family of vector fields

$$\mathcal{E}' = \left\{ X_f \mid f \in C^\infty(M)^G \right\} \tag{3.1}$$

spans the distribution $E$.  ♦

**3.2 Definition.**   *Let $(M, \{\cdot, \cdot\})$ be a Poisson manifold and $D \subset TM$ be a smooth generalized distribution on $M$. The distribution $D$ is called **Poisson** or **canonical**, if the condition $\mathbf{d}f|_D = \mathbf{d}g|_D = 0$ for $f, g \in C^\infty(M)$ implies that $\mathbf{d}\{f, g\}|_D = 0$.*

**3.3 Proposition.**   *The $G$–characteristic distribution $E$ is smooth, completely integrable, and Poisson. Its integral manifolds are given by the orbits of the pseudogroup $G_E$ of local diffeomorphisms generated by $\mathcal{E}$.*

**Proof.** The generalized distribution $E$ is smooth since it is spanned by all vector fields in $\mathcal{E}$. We will prove its complete integrability by using Theorem 2.3 which, at the same time, provides us with the characterization of the integral manifolds in terms of the $G_E$–orbits. So, let $m \in M$ and, for simplicity in the exposition, take $\mathcal{F}_T = F_T \in G_E$, with $F_t$ the flow of $X_f$, $f \in C^\infty(U)^G$, $U$ an open $G$–invariant neighborhood of $m$ (the general case in which $\mathcal{F}_T$ is the composition of a finite number of flows follows easily by attaching to what we are going to do a straightforward induction argument). Recall that the $G$–invariance of $f$ implies that $X_f$ and its flow $F_t$ are $G$–equivariant and consequently $\mathrm{Dom}(F_T)$ is a $G$–invariant open subset of $U$. The theorem follows if we are able to show that

$$T_m F_T(E(m)) = E(F_T(m)).$$

Let $X_g(m) \in E(m)$, with $g \in C^\infty(W)^G$, $W$ an open $G$–invariant subset of $\mathrm{Dom}(F_T)$. Since any Hamiltonian flow is always a Poisson map, by Proposition 2.6 we have that

$$T_m F_T(X_g(m)) = T_m F_T(X_{g \circ F_{-T} \circ F_T}(m)) = X_{g \circ F_{-T}}(F_T(m)) \in E(F_T(m)),$$

since $g \circ F_{-T} \in C^\infty(F_T(W))^G$ by the $G$–equivariance of $F_T$ and the $G$–invariance of $g$. This implies that $T_m F_T(E(m)) \subset E(F_T(m))$. Conversely, let

$$X_g(F_T(m)) \in E(F_T(m)).$$

Again, by using Proposition 2.6, $X_g(F_T(m)) = T_m F_T(X_{g \circ F_T}(m))$ which, by the equivariance of $F_T$, concludes the proof of the integrability of $E$ (for simplicity we omitted straightforward domain issues).

As to $E$ being canonical, let $f, g \in C^\infty(M)$ be such that

$$\mathbf{d}f|_E = \mathbf{d}g|_E = 0.$$

Consider $m \in M$ and let $h \in C^\infty(U)^G$ be arbitrary, with $U$ an open $G$–invariant neighborhood of $m$. Then,

$$\begin{aligned}
\mathbf{d}\{f, g\}(m) \cdot X_h(m) &= X_h\left[\{f, g\}\right](m) \\
&= \{\{f, g\}, h\}(m) \\
&= -\{\{h, f\}, g\}(m) - \{\{g, h\}, f\}(m) \\
&= \{X_h[f], g\}(m) - \{X_h[g], f\}(m) = 0,
\end{aligned}$$

as required.   ∎

**3.4 Definition.** *Let $(M, \{\cdot, \cdot\})$ be a Poisson manifold, $G$ be a Lie group acting canonically on $M$, and $E$ be the associated integrable $G$–characteristic distribution. Let $\mathcal{J} : M \to M/G_E$ be the canonical projection of $M$ onto the $G_E$–orbit space. We will call $\mathcal{J}$ the **optimal momentum map** of the canonical $G$–action on $(M, \{\cdot, \cdot\})$. We will refer to $M/E := M/G_E$ as the **momentum space**.*

A straightforward consequence of the previous definition is the following

**3.5 Theorem** (Optimal Noether Theorem). *Let $(M, \{\cdot, \cdot\})$ be a Poisson manifold, $G$ be a Lie group acting canonically on $M$, $E$ be the associated integrable $G$–characteristic distribution, and $\mathcal{J} : M \to M/G_E$ be the optimal momentum mapping. Then $\mathcal{J}$ is a constant of the motion for the dynamics generated by any $G$–invariant Hamiltonian $h$, that is,*

$$\mathcal{J} \circ F_t = \mathcal{J},$$

*where $F_t$ is the flow of $X_h$.*

**Remark.** By the very construction of $\mathcal{J}$, its level sets are the smallest immersed submanifolds (actually, we will see that under certain hypotheses they are embedded) respected by all $G$–equivariant Hamiltonian dynamics. This justifies the use of *optimal* in its denomination.  ♦

Notice that in contrast with the ordinary momentum map, $\mathcal{J}$ is always defined, which solves some of the problems of the traditional approach to the study of symmetries pointed out in the introduction. In particular, we have the following examples.

**Example. Conservation laws without momentum maps:** Consider the example presented in the introduction consisting of a canonical symmetry to which it is impossible to associate a globally defined momentum map. Let $M = S^1 \times S^1 = \mathbb{T}^2$ be the two torus with the symplectic form $\omega = d\theta_1 \wedge d\theta_2$. Consider $G = S^1$ acting on $M$ by

$$e^{i\phi} \cdot (e^{i\theta_1}, e^{i\theta_2}) := (e^{i(\phi+\theta_1)}, e^{i\theta_2}).$$

In order to compute the optimal momentum map $\mathcal{J}$, the first ingredient that we need is the $S^1$–characteristic distribution $E$. It is easy to see that in this case, every $S^1$–invariant function $f \in C^\infty(\mathbb{T}^2)^{S^1}$ can be written as

$$f(e^{i\theta_1}, e^{i\theta_2}) = g(e^{i\theta_2}),$$

with $g \in C^\infty(S^1)$. Its associated Hamiltonian vector field is given by $X_f = \frac{\partial g}{\partial \theta_2} \frac{\partial}{\partial \theta_1}$. Since $g$ is an arbitrary function on the circle, we can identify $M/E$ with the second circle $S^1$ in the torus $\mathbb{T}^2$. The optimal momentum map is therefore given by the expression:

$$\mathcal{J}: \begin{array}{ccc} \mathbb{T}^2 & \longrightarrow & S^1 \\ (e^{i\theta_1}, e^{i\theta_2}) & \longmapsto & e^{i\theta_2}. \end{array}$$

It is a remarkable fact that in this case the optimal momentum map is $S^1$–valued and moreover, it coincides with the Lie group valued momentum map introduced by McDuff [1988], Weitsman [1993], and Alekseev, Malkin, and Meinrenken [1997] that one would obtain by considering our example as a quasi–Hamiltonian $S^1$–space (see the prior references for an explanation of this term). ♦

**Example. The optimal momentum map of a Poisson non Hamiltonian action:** The previous example needed the introduction of group valued momentum maps in order to encode the conservation laws associated to the symmetries of the problem. We now give an example where even such momentum maps are not available. Nevertheless we will see that the optimal momentum map can carry out that job. Let $(\mathbb{R}^3, \{\cdot, \cdot\})$ be the Poisson manifold formed by the Euclidean three dimensional space $\mathbb{R}^3$ together with the Poisson structure induced by the Poisson tensor $B$ that in Euclidean coordinates takes the form:

$$B = \begin{pmatrix} 0 & 1 & 0 \\ -1 & 0 & 1 \\ 0 & -1 & 0 \end{pmatrix}.$$

With this Poisson bracket, the Hamiltonian vector field $X_f$ associated to any smooth function $f \in C^\infty(\mathbb{R}^3)$ is given by

$$X_f(x,y,z) = \frac{\partial f}{\partial y}\frac{\partial}{\partial x} + \left(\frac{\partial f}{\partial z} - \frac{\partial f}{\partial x}\right)\frac{\partial}{\partial y} - \frac{\partial f}{\partial y}\frac{\partial}{\partial z}. \qquad (3.2)$$

Consider the action of the additive group $(\mathbb{R}, +)$ on $\mathbb{R}^3$ given by $\lambda \cdot (x, y, z) := (x + \lambda, y, z)$, for any $\lambda \in \mathbb{R}$ and any $(x, y, z) \in \mathbb{R}^3$. In view of (3.2), it is clear that this action does not have a standard associated momentum map. Nevertheless, it is a Poisson action and therefore we can construct an optimal momentum map for it.

Indeed, notice first that the invariant functions $f \in C^\infty(M)^\mathbb{R}$ are all of the form $f(x, y, z) \equiv \bar{f}(y, z)$, with $\bar{f} \in C^\infty(\mathbb{R}^2)$ arbitrary. This implies that the $G_E$–orbits on $\mathbb{R}^3$ coincide with those of the $\mathbb{R}^2$–action on $\mathbb{R}^3$ given by $(\mu, \nu) \cdot (x, y, z) := (x + \mu, y + \nu, z - \mu)$, for any $(\mu, \nu) \in \mathbb{R}^2$ and any $(x, y, z) \in \mathbb{R}^3$. Therefore, $M/G_E$ can be identified with $\mathbb{R}$ and the associated optimal momentum map takes the form

$$\mathcal{J}: \begin{array}{ccc} \mathbb{R}^3 & \longrightarrow & \mathbb{R} \\ (x, y, z) & \longmapsto & x + z. \end{array}$$

It is easy to verify that the Hamiltonian flow associated to any invariant function $f(x,y,z) \equiv \bar{f}(y,z)$ preserves the level sets of $\mathcal{J}$; moreover, the function $\mathcal{J}$ is a Casimir of the Poisson manifold $(\mathbb{R}^3, \{\cdot,\cdot\})$. ♦

**Example. A canonical linear action:** Consider $\mathbb{C}^3$ with the symplectic form $\omega$ given by

$$\omega((z_1, z_2, z_3), (z_1', z_2', z_3')) = -\mathrm{Im}\,\langle (z_1, z_2, z_3), (z_1', z_2', z_3')\rangle.$$

Consider now the natural action of the Lie group SU(3) on $\mathbb{C}^3$ via matrix multiplication. This action is canonical and since it is linear, it is globally Hamiltonian. Moreover, given that the isotropy subgroup of any point in $\mathbb{C}^3$ with respect to this action has dimension at least three, the ordinary momentum map is always singular. We will see that its image can be naturally identified with the momentum space associated to the SU(3)–characteristic distribution.

Indeed, given that every SU(3)-invariant function in $\mathbb{C}^3$ is a function of

$$f(z_1, z_2, z_3) = \frac{1}{2}\left(|z_1|^2 + |z_2|^2 + |z_3|^2\right)$$

and the Hamiltonian flow of $X_f$ is given by

$$F_t(z_1, z_2, z_3) = (z_1 e^{-it}, z_2 e^{-it}, z_3 e^{-it}),$$

the orbit space $\mathbb{C}^3/G_E$ coincides with $\mathbb{C}^3/S^1$, where $S^1$ acts on $\mathbb{C}^3$ by

$$e^{i\phi} \cdot (z_1, z_2, z_3) = (e^{i\phi} z_1, e^{i\phi} z_2, e^{i\phi} z_3). \tag{3.3}$$

This quotient space can be identified with $(\mathbb{CP}(2) \times \mathbb{R}^+) \cup \{*\}$, where $\{*\}$ denotes a singleton or, said differently, with the cone $\overset{\circ}{C}(\mathbb{CP}(2))$ on $\mathbb{CP}(2)$. Indeed, if $\pi: \mathbb{C}^3 \to \mathbb{C}^3/S^1$ is the canonical projection and $\mathbf{z} = (z_1, z_2, z_3)$, then the mapping that assigns $\pi(z_1, z_2, z_3)$ to $([\mathbf{z}/\|\mathbf{z}\|], \|\mathbf{z}\|)$ if $\mathbf{z} \neq \mathbf{0}$, and to $*$ if $\mathbf{z} = \mathbf{0}$, provides the needed identification (the symbol $[\mathbf{z}/\|\mathbf{z}\|]$ denotes the element $\pi(\mathbf{z}/\|\mathbf{z}\|) \in \mathbb{CP}(2)$). We have the following expression for the optimal momentum map:

$$\begin{aligned}\mathcal{J}: \mathbb{C}^3 &\longrightarrow (\mathbb{CP}(2) \times \mathbb{R}^+) \cup \{*\} \\ \mathbf{z} &\longmapsto \begin{cases} \left(\left[\frac{\mathbf{z}}{\|\mathbf{z}\|}\right], \|\mathbf{z}\|\right) & \text{if } \mathbf{z} \neq \mathbf{0} \\ * & \text{if } \mathbf{z} = \mathbf{0}. \end{cases}\end{aligned}$$

♦

**Remark. The compact case and the Theory of Invariants:** The example that we just described lies in a very big class of systems for which the computation of the $G$–characteristic distribution $E$, and therefore of the optimal momentum map $\mathcal{J}$, is particularly simple. We are referring to canonical $G$–actions with $G$ a compact Lie group. It turns out that, according to a theorem of Gotay and Tuynman [1991], every canonical action of a compact Lie group on a symplectic manifold can be reduced to the study of a symplectic linear representation of $G$ on a certain finite dimensional vector space $V \simeq \mathbb{R}^{2n}$. Once we have reduced the problem to the linear representation of a compact Lie group, we have at our disposal the **Theory of Invariants**. For our purposes, the most interesting result in this theory is the Hilbert–Weyl Theorem (Weyl [1946]; Poènaru [1976]; Kempf [1987]) which guarantees that the algebra of $G$–invariant polynomials is finitely generated, that is, one can always find a finite number of $G$–invariant polynomials $\{\sigma_1, \ldots, \sigma_k\}$ such that every $G$–invariant polynomial $P \in \mathcal{P}(V)^G$ can be written as a polynomial function of them. More specifically, given $P \in \mathcal{P}(V)^G$, there is some $\widehat{P} \in \mathbb{R}[X_1, \ldots, X_k]$ such that $P = \widehat{P}(\sigma_1, \ldots, \sigma_k)$. Note that the generating family $\{\sigma_1, \ldots, \sigma_k\}$ can be chosen to be minimal. In that situation we say that $\{\sigma_1, \ldots, \sigma_k\}$ is a **Hilbert basis** of $\mathcal{P}(V)^G$. In applications, it is convenient to choose the Hilbert basis formed of homogeneous polynomials. Note also that the Hilbert basis is not necessarily free and that therefore there are in general relations between its elements.

The generalization of the Hilbert–Weyl Theorem to smooth functions has been carried out by Schwarz [1974], who proved that if $f$ is a germ of a function in $C^\infty(V)^G$ and $\{\sigma_1, \ldots, \sigma_k\}$ is a Hilbert basis of $\mathcal{P}(V)^G$, then there is a germ $\widehat{f} \in C^\infty(\mathbb{R}^k)$ such that $f = \widehat{f}(\sigma_1, \ldots, \sigma_k)$. Consequently, using (3.1) we can write the $G$–characteristic distribution in this case as

$$E = \text{span}\{X_{\sigma_1}, \ldots, X_{\sigma_k}\}. \qquad \blacklozenge$$

## 3.2 The optimal momentum map for proper globally Hamiltonian actions

In order to illustrate the content of $\mathcal{J}$, we now identify in the classical language the conservation laws induced by $\mathcal{J}$. In the following paragraphs, we assume that $(M, \omega)$ is a symplectic manifold and that $G$ is a Lie group acting properly on $M$ in a globally Hamiltonian fashion with associated momentum map $\mathbf{J} : M \to \mathfrak{g}^*$.

**3.6 Theorem.**   *Let $(M, \omega)$ be a symplectic manifold and $G$ be a Lie group acting properly on $M$ in a globally Hamiltonian fashion with associated momentum map $\mathbf{J} : M \to \mathfrak{g}^*$. If $E$ is the $G$–characteristic distribution, then for any $m \in M$:*

$$E(m) = \ker T_m \mathbf{J} \cap T_m M_{G_m}. \tag{3.4}$$

*Moreover, the $G_E$–orbit of the point $m$, and therefore the level set $\mathcal{J}^{-1}(\rho)$, $\rho \in M/E$, of $\mathcal{J}$ containing the point $m \in M$, is the connected component* $(\mathbf{J}^{-1}(\mu) \cap M_H)_{c.c.m}$ *containing m of the embedded submanifold* $\mathbf{J}^{-1}(\mu) \cap M_H$, *that is:*

$$\mathcal{J}^{-1}(\rho) = (\mathbf{J}^{-1}(\mu) \cap M_H)_{c.c.m}$$

*where* $\mu = \mathbf{J}(m) \in \mathfrak{g}^*$ *and* $H := G_m$ *is the isotropy subgroup of* $m \in M$.

**Proof.** Expression (3.4) is a consequence of the following chain of equalities:

$$\begin{aligned}
E(m) &= \mathrm{span}\{X_f(m) | f \in C^\infty(M)^G\} \\
&= B^\sharp(m)\left(\mathrm{span}\{\mathbf{d}f(m) | f \in C^\infty(M)^G\}\right) \\
&= B^\sharp(m)\left(\left((T_m(G \cdot m))^\circ\right)^{G_m}\right) &\text{(by Proposition 2.14)} \\
&= \left(B^\sharp(m)\left((T_m(G \cdot m))^\circ\right)\right)^{G_m} &\text{(by Proposition 2.17)} \\
&= ((T_m(G \cdot m))^\omega)^{G_m} &\text{(by Proposition 2.4)} \\
&= (\ker T_m\mathbf{J})^{G_m} = \ker T_m\mathbf{J} \cap T_m M_{G_m}.
\end{aligned}$$

We now prove the claim in the statement about the integral manifolds of $E$. Let again $m \in M$ be such that $\mu = \mathbf{J}(m) \in \mathfrak{g}^*$, $\mathcal{J}(m) = \rho$, and $H := G_m \subset G$ is its isotropy subgroup. As we said in the introduction, the subset $M_H \subset M$ is a symplectic submanifold of $M$. Moreover, it is easy to see Otto [1987]; Ortega [1998]; Ortega and Ratiu [2002] that the restriction $\mathbf{J}|_{M_H}$ of $\mathbf{J}$ to $M_H$ is a constant rank mapping and hence, by the Fibration Theorem [Abraham, Marsden, and Ratiu, 1988, Theorem 3.5.18], $\mathbf{J}|_{M_H}^{-1}(\mu) = \mathbf{J}^{-1}(\mu) \cap M_H$ is a submanifold of $M_H$ and consequently of $M$, which contains the point $m$. Given that for any point $z \in \mathbf{J}^{-1}(\mu) \cap M_H$ we have that $T_z\left(\mathbf{J}^{-1}(\mu) \cap M_H\right) = T_z\left(\mathbf{J}|_{M_H}^{-1}(\mu)\right) = \ker T_z\mathbf{J} \cap T_z M_H = E(z)$, we can conclude that the connected component $(\mathbf{J}^{-1}(\mu) \cap M_H)_{c.c.m}$ containing the point $m$ of the submanifold $\mathbf{J}^{-1}(\mu) \cap M_H$ is an integral manifold of $E$ containing $m$. At the same time, the characterization of the level sets of $\mathcal{J}$ as $G_E$–orbits implies, via the standard Noether's Theorem and the principle of conservation of the isotropy that $\mathcal{J}^{-1}(\rho) = G_E \cdot m \subset (\mathbf{J}^{-1}(\mu) \cap M_H)_{c.c.m}$. The result follows from the uniqueness of the maximal integral manifolds of a generalized distribution [Libermann and Marle, 1987, Theorem 2.3, Appendix 3]. ∎

**Remark.** The previous theorem justifies again the use of the adjective *optimal* in the denomination of $\mathcal{J}$ since it proves that in the proper globally Hamiltonian case, its level sets coincide with the smallest invariant subsets of $M$ under $G$–equivariant dynamics. In other words, the optimal momentum map $\mathcal{J}$ is capable of implementing in one shot both the classical Noether Theorem, as well as the law of conservation of the isotropy. ♦

The use of the name *momentum* for $\mathcal{J}$ is reasonable since, as it follows from Theorem 3.6, there are cases in which $\mathcal{J}$ and $\mathbf{J}$ are basically the same map. For example, suppose that we are in the hypotheses of Theorem 3.6 and, additionally, we assume the $G$–action to be free (there are no singularities) and the level sets of $\mathbf{J}$ connected. In this situation, the map $\mathcal{J}(m) = \rho \mapsto \mathbf{J}(m) = \mu$ is a bijection $\varphi$ between $M/G_E$ and $\mathfrak{g}_{\mathbf{J}}^* := \mathbf{J}(M)$. Indeed, it is well defined since if we take another $m' \in M$ such that $\mathcal{J}(m') = \rho$, then $m$ and $m'$ are in the same level set of $\mathcal{J}$ that in our hypotheses, by Theorem 3.6, are the (connected) level sets of $\mathbf{J}$ and hence $\varphi(\mathcal{J}(m')) = \mathbf{J}(m') = \mu = \mathbf{J}(m) = \varphi(\mathcal{J}(m))$. The map $\varphi$ is onto by construction and one–to–one by the connectedness of the level sets of $\mathbf{J}$. Note that in this case, the commutative diagram provides an identification between $\mathcal{J}$ and $\mathbf{J}$. This diagram is a corollary of the universality properties of the optimal momentum map that we study in the following subsection.

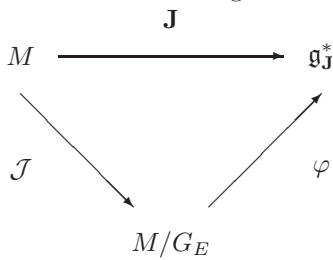

## 3.3 Universality properties of the optimal momentum map

In this section we will show that the optimal momentum map is universal in the sense of Category Theory, that is, any other momentum map that we may define is going to factor through $\mathcal{J}$. Before making this statement more explicit we will need to introduce a few properties of the orbit space $M/G_E$ and the optimal momentum map. We start with the following.

**3.7 Definition.** *A pair $(X, C^\infty(X))$, where $X$ is a topological space and $C^\infty(X) \subset C^0(X)$ is a subset of continuous functions in $X$, is called a **variety** with **smooth functions** $C^\infty(X)$. If $Y \subset X$ is a subset of $X$, the pair $(Y, C^\infty(Y))$ is said to be a **subvariety** of $(X, C^\infty(X))$, if $Y$ is a topological space endowed with the relative topology defined by that of $X$ and*

$$C^\infty(Y) = \{f \in C^0(Y) \mid f = F|_Y \text{ for some } F \in C^\infty(X)\}.$$

*Sometimes $C^\infty(Y)$ is called the set of **Whitney smooth functions** on $Y$*

with respect to $X$. A map $\varphi : X \to Z$ between two varieties is said to be smooth when it is continuous and $\varphi^* C^\infty(Z) \subset C^\infty(X)$.

In our discussion we are interested in the varieties constructed using generalized integrable distributions $D$ on the manifold $M$. If we denote by $M/D$ the space formed by the integral manifolds of $M$, the pair $(M/D, C^\infty(M/D))$ is a variety whose set of smooth functions $C^\infty(M/D)$ is defined by the requirement that the canonical projection $\pi : M \to M/D$ is a smooth map, that is,

$$C^\infty(M/D) := \{f \in C^0(M/D) \mid f \circ \pi \in C^\infty(M)\}.$$

Note that if $(M, \{\cdot, \cdot\})$ is a Poisson manifold and $D \subset TM$ is a Poisson integrable distribution, the pair $(C^\infty(M/D), \{\cdot, \cdot\}_{M/D})$ is a well–defined Poisson algebra [Ortega and Ratiu, 1998, Theorem 2.12], where the bracket $\{\cdot, \cdot\}_{M/D}$ is given by

$$\{f, g\}_{M/D}(\pi(m)) = \{f \circ \pi, g \circ \pi\}(m), \tag{3.5}$$

for every $m \in M$. In the particular case in which $E$ is the $G$–characteristic distribution, Proposition 3.3 guarantees that $(C^\infty(M/E), \{\cdot, \cdot\}_{M/E}) := (C^\infty(M/G_E), \{\cdot, \cdot\}_{M/G_E})$ is a well–defined Poisson algebra. Moreover, the construction of the bracket (3.5) implies that the optimal momentum map $\mathcal{J}$ is a smooth Poisson morphism.

**3.8 Proposition.**   Let $(M, \{\cdot, \cdot\})$ be a Poisson manifold, $G$ be a Lie group acting canonically on it, and $E$ and $\mathcal{J}$ be the associated $G$–characteristic distribution and optimal momentum map, respectively. Let $m \in M$ be arbitrary such that $\mathcal{J}(m) = \rho \in M/G_E$. Then, for any $g \in G$, the map $\Phi_g(\rho) = \mathcal{J}(g \cdot m) \in M/G_E$ defines a smooth Poisson $G$–action on $M/G_E$ with respect to which $\mathcal{J}$ is $G$–equivariant.

**Proof.** We just have to check that $\Phi$ is well defined since if that is case, the equivariance of $\mathcal{J}$ will follow by construction. Let $m, m' \in M$ be such that $\mathcal{J}(m) = \mathcal{J}(m') = \rho$. This implies that $m$ and $m'$ live in the same integral manifold of $E$, that is, in the same $G_E$–orbit. Hence, there exists $\mathcal{F}_T \in G_E$ such that $\mathcal{F}_T(m) = m'$. Since $\mathcal{F}_T$ is the composition of a finite number of $G$–equivariant Hamiltonian flows associated to $G$–invariant Hamiltonians, it is $G$–equivariant and therefore

$$\mathcal{J}(g \cdot m') = \mathcal{J}(g \cdot \mathcal{F}_T(m)) = \mathcal{J}(\mathcal{F}_T(g \cdot m)) = \mathcal{J}(g \cdot m) = \Phi_g(\rho),$$

as required. The smoothness and the Poisson character of $\Phi$ follow from the fact that this action is the projection onto $M/E$ of the smooth Poisson action on the manifold $M$, via the smooth optimal momentum map.  ∎

**Example.**   We look at the Poisson and $G$–structures of the spaces $M/E$ found in the examples 3.1 and 3.1.

- $S^1$ **acting on** $\mathbb{T}^2$: In this case $M/E$ coincides with $S^1$ and the optimal momentum map is given by $\mathcal{J}(e^{i\theta_1}, e^{i\theta_2}) = e^{i\theta_2}$. The group $S^1$ acts trivially on $M/E \simeq S^1$ and the Poisson structure $\{\cdot, \cdot\}_{M/E}$ is trivial. Indeed, let $f, g \in C^\infty(S^1)$ and $(e^{i\theta_1}, e^{i\theta_2}) \in \mathbb{T}^2$ be arbitrary. Then

$$\begin{aligned}\{f,g\}_{M/E}(e^{i\theta_2}) &= \{f,g\}_{M/E}(\mathcal{J}(e^{i\theta_1}, e^{i\theta_2})) \\ &= \{f \circ \mathcal{J}, g \circ \mathcal{J}\}_{\mathbb{T}^2}(e^{i\theta_1}, e^{i\theta_2}) \\ &= d\theta_1 \wedge d\theta_2 \left(\frac{\partial f}{\partial \theta_2}\frac{\partial}{\partial \theta_1}, \frac{\partial g}{\partial \theta_2}\frac{\partial}{\partial \theta_1}\right) = 0.\end{aligned}$$

- $\mathrm{SU}(3)$ **acting on** $\mathbb{C}^3$: As we saw, $M/E \simeq (\mathbb{CP}(2) \times \mathbb{R}^+) \cup \{*\}$. The Lie group $\mathrm{SU}(3)$ acts on this set by leaving the point $*$ fixed and by

$$A \cdot \left(\left[\frac{\mathbf{z}}{\|\mathbf{z}\|}\right], \|\mathbf{z}\|\right) = \left(\left[\frac{A\mathbf{z}}{\|A\mathbf{z}\|}\right], \|A\mathbf{z}\|\right),$$

when $A \in \mathrm{SU}(3)$ and $\mathbf{z} \neq 0$. As to the Poisson structure on $M/E$, it can be easily shown that for any $f, g \in C^\infty(M/E)$ and any $([z], r) \in M/E \simeq (\mathbb{CP}(2) \times \mathbb{R}^+) \cup \{*\}$

$$\{f,g\}_{M/E}([z], r) = \omega_{\mathbb{CP}(2)}([z])\left(X_{f_r}([z]), X_{g_r}([z])\right),$$

and

$$\{f,g\}_{M/E}(*) = 0,$$

where $\omega_{\mathbb{CP}(2)}$ in the natural symplectic structure on $\mathbb{CP}(2)$ (coming from considering it as one of the regular symplectic reduced spaces of the $S^1$–action (3.3) on $\mathbb{C}^3$) and $f_r, g_r \in C^\infty(\mathbb{CP}(2))$ are defined by $f_r([z]) := f([z], r)$ and $g_r([z]) := g([z], r)$, for any $[z] \in \mathbb{CP}(2)$.  ♦

In order to illustrate the universality properties of the optimal momentum map we introduce the category of Hamiltonian symmetric systems with a momentum map.

**3.9 Definition.** *Let $(M, \{\cdot, \cdot\})$ be a Poisson manifold, $(P, C^\infty(P), \{\cdot, \cdot\}_P)$ be a Poisson variety, and $G$ be a Lie group acting canonically on $M$ and $P$. Let $\mathbf{K} : M \to P$ be a smooth $G$–equivariant Poisson map. We say that $\mathbf{K}$ is a momentum map for the $G$–action on $M$ if the Hamiltonian flows associated to $G$–invariant smooth functions leave invariant the level sets of $\mathbf{K}$. In that situation we say that $(M, \{\cdot, \cdot\}, G, \mathbf{K} : M \to P)$ is a **Hamiltonian $G$–space with momentum map $\mathbf{K}$**.*

**3.10 Theorem** (Universality of the optimal momentum map). *The optimal momentum map is a universal object in the category of Hamiltonian symmetric systems with a momentum map. More specifically, if*

$$(M, \{\cdot,\cdot\}, G, \mathbf{K} : M \to P)$$

*is any Hamiltonian G–space with momentum map $\mathbf{K}$ and $\mathcal{J} : M \to M/E$ is the optimal momentum map defined using the canonical G–action on $M$, then there exists a unique G–equivariant Poisson morphism $\varphi : M/E \to P$ such that the diagram commutes.*

$$\begin{array}{ccc} M & \xrightarrow{\mathbf{K}} & P \\ & \searrow{\mathcal{J}} \quad \nearrow{\varphi} & \\ & M/E & \end{array}$$

**Proof.** The function $\varphi$ is given, for any $\rho = \mathcal{J}(m) \in M/E$, by the expression

$$\varphi(\rho) := \mathbf{K}(m).$$

The map $\varphi$ is well defined since if $m' \in \mathcal{J}^{-1}(\rho)$, then there exists a finite composition $\mathcal{F}_T$ of flows associated to $G$–invariant Hamiltonians such that $m' = \mathcal{F}_T(m)$. Given that $\mathbf{K}$ is a momentum map we have that

$$\mathbf{K}(m') = \mathbf{K}(\mathcal{F}_T(m)) = \mathbf{K}(m) = \varphi(\rho).$$

The smoothness of $\varphi$, as well as its $G$–equivariance, and Poisson character are a simple diagram chasing exercise. The uniqueness is guaranteed by the fact that the diagram (3.10) commutes and by the surjectivity of $\mathcal{J}$. ∎

## 4  Optimal reduction

As we already said in the introduction, the most efficient way to profit from the conservation laws encoded in the symmetries of a globally Hamiltonian system is carrying out the so called symplectic or Marsden–Weinstein reduction (Marsden and Weinstein [1974]), which we briefly review. Let $(M, \omega, h, G, \mathbf{J} : M \to \mathfrak{g}^*)$ be a globally Hamiltonian symmetric system where we will assume that the $G$–action on $M$ is free and proper. Let $\mu \in \mathfrak{g}^*_\mathbf{J}$ be an arbitrary element in the image of the momentum map $\mathbf{J}$. The Marsden–Weinstein reduction theorem says that the quotient $\mathbf{J}^{-1}(\mu)/G_\mu$ is a symplectic manifold with symplectic form $\omega_\mu$ uniquely determined by the equality $\pi_\mu^* \omega_\mu = i_\mu^* \omega$, where $i_\mu : \mathbf{J}^{-1}(\mu) \hookrightarrow M$ and $\pi_\mu : \mathbf{J}^{-1}(\mu) \to \mathbf{J}^{-1}(\mu)/\pi_\mu$ are the natural injection and projection respectively. The dynamics induced by $X_h$ projects naturally onto the reduced space $\mathbf{J}^{-1}(\mu)/G_\mu$.

In this section we see how the new formulation in terms of the optimal momentum map allows us to mimic this procedure, creating the possibility to reduce symmetric systems in all the situations in which $\mathcal{J}$ is defined and freeing us from the strong restrictions posed by the classical formulations of the reduction theorems.

## 4.1 Reduction lemmas

The first ingredient needed in the reduction of a symmetric system are the level sets of the associated momentum map. Since by construction the level sets of the optimal momentum map $\mathcal{J}$ are the integral manifolds of a smooth distribution, they are always smooth immersed submanifolds of $M$. Moreover, in Theorem 3.6 we saw that if $\mathcal{J}$ is associated to a proper globally Hamiltonian action, its level sets are actually embedded submanifolds. In the next result we show that this is also the case under much weaker hypotheses. We start with the following straightforward lemma.

**4.1 Lemma.** *Let $(M, \omega)$ be a symplectic manifold and $G$ be a Lie group acting properly and canonically on $M$. Let $E$ be the $G$–characteristic distribution with optimal momentum map $\mathcal{J} : M \to M/G_E$. Then, for any $\rho \in M/G_E$, the set $\mathcal{J}^{-1}(\rho) \subset M$ is included in the connected component of some isotropy type manifold $M_H$, with $H$ the isotropy subgroup of any $m \in \mathcal{J}^{-1}(\rho)$.*

**Proof.** It is a straightforward consequence of the equality $\mathcal{J}^{-1}(\rho) = G_E \cdot m$, for any $m \in \mathcal{J}^{-1}(\rho)$. ∎

**4.2 Definition.** *In the hypotheses of the previous lemma, we say that an element $\rho \in M/G_E$ satisfies the **closedness hypothesis** if $\mathcal{J}^{-1}(\rho)$ is closed as a subset of the isotropy type submanifold $M_H$ in which it is sitting.*

**Example.** The closedness hypothesis is always satisfied in the presence of globally Hamiltonian actions. Also, suppose that $\rho \in M/E$ is such that $\mathcal{J}^{-1}(\rho) \subset M_H$. Let $N(H)$ be the normalizer in $G$ of $H$. The canonical $G$–action on $M$ induces a natural canonical action of the group $N(H)/H$ on the symplectic manifold $M_H$. Let $(N(H)/H)^\rho$ be the subgroup of $N(H)/H$ that leaves invariant the connected component $M_H^\rho$ of $M_H$ in which $\mathcal{J}^{-1}(\rho)$ is sitting. If the $(N(H)/H)^\rho$–action on $M_H^\rho$ has a globally defined momentum map associated, then $\rho$ satisfies the closedness hypothesis. ♦

**4.3 Proposition.** *Let $(M, \omega)$ be a symplectic manifold and $G$ be a Lie group acting properly and canonically on $M$. Let $E$ be the $G$–characteristic distribution with optimal momentum map $\mathcal{J} : M \to M/G_E$. Let $\rho$ be an element in $M/G_E$ that satisfies the closedness hypothesis. If $M_H$ is the isotropy type submanifold in which the level set $\mathcal{J}^{-1}(\rho)$ is included by Lemma 4.1, then $\mathcal{J}^{-1}(\rho)$ is a closed embedded submanifold of $M_H$ and therefore an embedded submanifold of $M$. As a consequence, if $M_H$ is closed in $M$ then $\mathcal{J}^{-1}(\rho)$ is a closed embedded submanifold of $M$.*

**Proof.** Let $M_H^\rho$ be the connected component of $M_H$ containing $\mathcal{J}^{-1}(\rho)$. Thus, the claim of the proposition will follow if we are able to show that $\mathcal{J}^{-1}(\rho)$ is a closed embedded submanifold of $M_H^\rho$.

Firstly, note that $M_H^\rho$ inherits from the canonical $G$–action on $M$ a free and canonical $L^\rho := (N(H)/H)^\rho$–action, where $N(H)$ denotes the normalizer of $H$ in $G$ and $(N(H)/H)^\rho$ is the closed subgroup of $N(H)/H$ that leaves $M_H^\rho$ invariant. In our subsequent discussion we will assume, in order to simplify the exposition, that $M_H^\rho = M_H$ and $(N(H)/H)^\rho = N(H)/H$. Let $E_L$ be the $L$-characteristic distribution associated to the canonical $L$-action on $M_H$ and $\mathcal{J}_L : M_H \to M_H/E_L$ be the associated optimal momentum map. Let $m \in \mathcal{J}^{-1}(\rho)$. In the sequel we will prove that if $\mathcal{J}_L(m) = \sigma$ then, $\mathcal{J}^{-1}(\rho) = \mathcal{J}_L^{-1}(\sigma)$. We first need the following lemma.

**4.4 Lemma.** *Let $G$ be a Lie group acting properly on the manifold $M$. Let $m \in M$ be an arbitrary point such that $H := G_m$. Then every function $f \in C^\infty(M_H)^{N(H)} = C^\infty(M_H)^L$ admits a local extension at $m$ to $C^\infty(M)^G$, that is, there exists a $G$–invariant neighborhood $U$ of $m$ in $M$ and a $G$-invariant function $F \in C^\infty(M)^G$, such that $F|_{U \cap M_H} = f|_{U \cap M_H}$.*

**Proof.** Since the claim of the lemma is local we can make use of the slices introduced in Section 2. Let $V$ be an open $G$–invariant neighborhood of the orbit $G \cdot m$ that is modeled by the tube $G \times_H A_r$. It is easy to see that

$$V_H \simeq N(H) \times_H A_r^H.$$

Let now $g : A_r \to \mathbb{R}$ be the smooth function defined by $g(v) := f([e, v])$ for any $v \in A_r^H$. Using a bump function similar to the one that was used in the proof of Proposition 2.13 we can construct another function $g_1 \in C^\infty(A_r^H)$ such that $g_1|_{A_{r/2}^H} = g|_{A_{r/2}^H}$ and $g_1|_{A_r^H \setminus A_{3r/4}^H} = 0$. Due to the compactness of $H$, the vector space $A$ can be decomposed as the direct sum $A = A^H \oplus W$ of two $H$–invariant subspaces $A^H$ and $W$. Define $g_2 \in C^\infty(A_r)^H$ by $g_2(v + w) = g_1(v)$, for any $v \in A_r^H$ and $w \in W$. We now let $g_3 \in C^\infty(V)^H$ be given by $g_3([h, v]) = g_2(v)$, for any $[h, v] \in V \simeq G \times_H A_r$. Finally, let $F \in C^\infty(M)^G$ be the function given by

$$F(z) = \begin{cases} g_3(z) & \text{if} \quad z \in V \\ 0 & \text{if} \quad z \notin V. \end{cases}$$

The lemma follows by taking the function $F$ above and $U$ as the open $G$–invariant set modeled by $G \times_H A_{r/2}$. ▼

As a corollary to this lemma we have that $E|_{M_H} = E_L$ and, consequently $\mathcal{J}^{-1}(\rho) = \mathcal{J}_L^{-1}(\sigma)$, as required.

We now show that the distribution $E_L$ has constant rank. Indeed, for any $z \in M_H$ we have that

$$\begin{aligned} E_L(z) &= \{X_f(z) \mid f \in C^\infty(M_H)^L\} \\ &= B_{M_H}^\sharp(z)\left(\{\mathbf{d}f(z) \mid f \in C^\infty(M_H)^L\}\right) \\ &= B_{M_H}^\sharp(z)\left(T_z(L \cdot z)^\circ\right) \qquad \text{(by Proposition 2.14)} \\ &= (T_z(L \cdot z))^{\omega|_{M_H}}. \end{aligned}$$

In particular, this equality shows that $E_L$ is an integrable distribution of constant rank equal to $\dim M_H - \dim \mathfrak{l}^*$. In the previous chain of equalities we denoted by $B_{M_H}$ the Poisson tensor associated to the symplectic form $\omega|_{M_H}$ on $M_H$.

The proof is concluded by recalling the closedness hypothesis on $\rho$ and a general fact about constant rank smooth foliations (see for instance Theorem 5 in page 51 of Camacho and Lins Neto [1985]) which states that the closed integral leaves of an integrable distribution of constant rank are always embedded submanifolds. This fact proves that $\mathcal{J}_L^{-1}(\sigma) = \mathcal{J}^{-1}(\rho)$ is an embedded submanifold of $M_H$, and thereby of $M$, as required. ■

**4.5 Proposition.** *In the hypotheses of the previous proposition, for any $\rho \in M/G_E$ satisfying the closedness hypothesis, the isotropy subgroup $G_\rho \subset G$ of $\rho$ with respect to the $G$–action on $M/G_E$ defined in Proposition 3.8, is a closed Lie subgroup of $G$. Moreover, for any $m \in \mathcal{J}^{-1}(\rho)$*

$$T_m(G_\rho \cdot m) = T_m(\mathcal{J}^{-1}(\rho)) \cap T_m(G \cdot m). \tag{4.1}$$

**Proof.** In order to show that $G_\rho$ is a Lie subgroup of $G$ it suffices to show (see for instance [Warner, 1983, Theorem 3.42]) that $G_\rho$ is closed in $G$. Let $\{g_n\} \subset G_\rho$ be an arbitrary convergent sequence in $G_\rho$ with limit $g \in G$. The closedness of $G_\rho$ will be guaranteed if we show that the limit $g$ actually belongs to $G_\rho$. Let $m \in M$ be such that $\mathcal{J}(m) = \rho$ and $H := G_m$. The condition $\{g_n\} \subset G_\rho$ implies that for any given $n \in \mathbb{N}^*$, there exists an element $\mathcal{F}_{T_n}^n \in G_E$ such that $g_n \cdot m = \mathcal{F}_{T_n}^n(m)$. Consequently, the $G$-equivariance of the elements $\mathcal{F}_{T_n}^n$ implies that the isotropy subgroups $G_{g_n \cdot m}$ of $g_n \cdot m$ satisfy that

$$g_n H g_n^{-1} = G_{g_n \cdot m} = G_{\mathcal{F}_{T_n}^n(m)} = G_m = H,$$

and hence the sequence $\{g_n\} \subset N(H)$. Since the normalizer $N(H)$ is closed in $G$, the limit $g$ belongs to $N(H)$ and hence the element $g \cdot m$ is sitting in the same connected component of $M_H$ in which the sequence $\{g_n \cdot m\}$ lives. Consequently, the element $g \cdot m$ lies in the closure in $M_H$ of $\mathcal{J}^{-1}(\rho)$. The closedness hypothesis on $\rho$ guarantees that $g \cdot m \in \mathcal{J}^{-1}(\rho) \subset M_H$. In particular, this implies that $g \cdot \rho = g \cdot \mathcal{J}(m) = \mathcal{J}(g \cdot m) = \rho$ or, equivalently, $g \in G_\rho$, as required.

We now prove equality (4.1). The inclusion $T_m(G_\rho \cdot m) \subset T_m(\mathcal{J}^{-1}(\rho)) \cap T_m(G \cdot m)$ is straightforward since the orbit $G_\rho \cdot m$ is included in both $\mathcal{J}^{-1}(\rho)$ and $G \cdot m$. Conversely, let

$$X_f(m) = \xi_M(m) \in T_m(\mathcal{J}^{-1}(\rho)) \cap T_m(G \cdot m), \tag{4.2}$$

with $f \in C^\infty(M)^G$ and $\xi \in \mathfrak{g}$. Recall that since the $G$–action on $M$ is canonical, the vector field $\xi_M \in \mathfrak{X}(M)$ is locally Hamiltonian and therefore there is a smooth function, say $\mathbf{J}^\xi \in C^\infty(U)$, locally defined in an open neighborhood $U$ of $m$ in $M$, such that $X_{\mathbf{J}^\xi} = \xi_M$. Notice that for any $z \in U$,

$$\{f, \mathbf{J}^\xi\}(z) = X_{\mathbf{J}^\xi}[f](z) = \mathbf{d}f(z) \cdot X_{\mathbf{J}^\xi}(z) = \mathbf{d}f(z) \cdot \xi_M(z) = 0,$$

by the $G$–invariance of the function $f$. Consequently, at any point in $U$, the Lie bracket $[X_f, X_{\mathbf{J}^\xi}] = X_{\{\mathbf{J}^\xi, f\}} = 0$, and hence, if $F_t$ is the flow of $X_f$ and $G_t$ is the flow of $X_{\mathbf{J}^\xi}$ (more explicitly $G_t(z) = \exp t\xi \cdot z$ for any $z \in U$), then $F_t \circ G_s = G_s \circ F_t$ (see for instance [ Abraham, Marsden, and Ratiu, 1988, Proposition 4.2.27]). By one of the Trotter product formulas (see Trotter [1958] or [ Abraham, Marsden, and Ratiu, 1988, Corollary 4.1.27]), the flow $H_t$ of $X_{f-\mathbf{J}^\xi} = X_f - X_{\mathbf{J}^\xi}$ is given by

$$\begin{aligned}H_t(z) &= \lim_{n\to\infty} \left(F_{t/n} \circ G_{-t/n}\right)^n(z) \\ &= \lim_{n\to\infty} \left(F^n_{t/n} \circ G^n_{-t/n}\right)(z) \\ &= (F_t \circ G_{-t})(z) = F_t(\exp -t\xi \cdot z),\end{aligned}$$

for any $z \in U$. Note that by (4.2), the point $m \in M$ is an equilibrium of $X_{f-\mathbf{J}^\xi} = X_f - X_{\mathbf{J}^\xi}$, hence $F_t(\exp -t\xi \cdot m) = m$ or, analogously $\exp t\xi \cdot m = F_t(m)$. Applying $\mathcal{J}$ on both sides of this equality, taking into account that $F_t$ is the flow of a $G$–invariant Hamiltonian vector field, it follows that $\exp t\xi \cdot \rho = \rho$, and hence $\xi \in \mathfrak{g}_\rho$. Thus $\xi_M(m) \in T_m(G_\rho \cdot m)$, as required. ∎

## 4.2   The optimal reduction method

We continue the study of the ingredients needed for reduction with the following.

**4.6 Proposition.**   *Let $(M, \omega)$ be a symplectic manifold and $G$ be a Lie group acting properly and canonically on $M$. Let $E$ be the $G$–characteristic distribution with optimal momentum map $\mathcal{J} : M \to M/G_E$. Then, for* any $\rho \in M/G_E$ satisfying the closedness hypothesis, the isotropy subgroup $G_\rho$ acts on the submanifold $\mathcal{J}^{-1}(\rho)$, and the corresponding orbit space $M_\rho := \mathcal{J}^{-1}(\rho)/G_\rho$ is a regular quotient manifold, that is, it can be endowed with the unique smooth structure that makes the canonical projection $\pi_\rho : \mathcal{J}^{-1}(\rho) \to \mathcal{J}^{-1}(\rho)/G_\rho$ a submersion. We will call $M_\rho = \mathcal{J}^{-1}(\rho)/G_\rho$ endowed with this smooth structure the **reduced phase space**.

**Proof.** Let $m \in \mathcal{J}^{-1}(\rho) \subset M$ be such that $H := G_m$. Recall that in the proof of Proposition 4.3 we showed that $\mathcal{J}^{-1}(\rho) \subset M_H$. This implies that the isotropies of all the elements of $\mathcal{J}^{-1}(\rho)$ under the $G_\rho$ action are identical and equal to $H \cap G_\rho$. A classical result (see for instance Exercise 4.1M in Abraham and Marsden [1978]) guarantees that in such situation the quotient $\mathcal{J}^{-1}(\rho)/G_\rho$ is a regular manifold, and the claim follows. ∎

We are now in position to state the main result of this section.

**4.7 Theorem** (Optimal Reduction). *Let $(M, \omega)$ be a symplectic manifold and $G$ be a Lie group acting properly and canonically on $M$. Let $E$ be the $G$–characteristic distribution with optimal momentum map $\mathcal{J} : M \to M/G_E$. Then, for any $\rho \in M/G_E$ satisfying the closedness hypothesis, the reduced space $M_\rho = \mathcal{J}^{-1}(\rho)/G_\rho$ has a unique symplectic structure $\omega_\rho$ characterized by*

$$\pi_\rho^* \omega_\rho = i_\rho^* \omega, \qquad (4.3)$$

*where $\pi_\rho : \mathcal{J}^{-1}(\rho) \to M_\rho$ is the canonical projection and $i_\rho : \mathcal{J}^{-1}(\rho) \to M$ is the inclusion.*

**Proof.** Let $[z]_\rho = \pi_\rho(z)$ be an arbitrary element of $M_\rho$. Since by Proposition 4.6 the projection $\pi_\rho$ is a surjective submersion, every vector $[v]_\rho \in T_{[z]_\rho} M_\rho$ can be written as $[v]_\rho = T_z \pi_\rho \cdot v$, with $v \in T_z(\mathcal{J}^{-1}(\rho))$. Taking also $[w]_\rho = T_z \pi_\rho \cdot w \in T_{[z]_\rho} M_\rho$ arbitrary, we define

$$\omega_\rho([z]_\rho)(T_z \pi_\rho \cdot v, T_z \pi_\rho \cdot w) := \omega(z)(v, w). \qquad (4.4)$$

In order to verify that this is a good definition we have to verify that it is independent of the representative $z \in \mathcal{J}^{-1}(\rho)$ that defines $[z]_\rho$ and of the vectors $v, w \in T_z(\mathcal{J}^{-1}(\rho))$ that define $[v]_\rho, [w]_\rho \in T_{[z]_\rho} M_\rho$, respectively. So, let $z' \in \mathcal{J}^{-1}(\rho)$ and $v', w' \in T_{z'}(\mathcal{J}^{-1}(\rho))$ be such that $[z]_\rho = [z']_\rho$ and $[v]_\rho = [v']_\rho, [w]_\rho = [w']_\rho$. Let $g \in G_\rho$ be such that $z' = g \cdot z$. Note that since $\pi_\rho = \pi_\rho \circ \Phi_g$ implies $T_z \pi_\rho = T_{z'} \pi_\rho \circ T_z \Phi_g$, the relation $[v]_\rho = [v']_\rho$ can be written as $T_{z'} \pi_\rho(v') = T_{z'} \pi_\rho(T_z \Phi_g(v))$. Consequently,

$$v' - T_z \Phi_g(v) \in \ker T_{z'} \pi_\rho = T_{z'}(G_\rho \cdot z')$$

and therefore there are elements $\xi^1, \xi^2 \in \mathfrak{g}_\rho$ such that

$$v' = T_z \Phi_g(v) + \xi^1_{\mathcal{J}^{-1}(\rho)}(z') \quad \text{and} \quad w' = T_z \Phi_g(w) + \xi^2_{\mathcal{J}^{-1}(\rho)}(z').$$

We now prove that $\omega_\rho([z']_\rho)([v']_\rho, [w']_\rho) = \omega_\rho([z]_\rho)([v]_\rho, [w]_\rho)$:

$$\begin{aligned}\omega_\rho([z']_\rho)([v']_\rho, [w']_\rho) &= \omega(\Phi_g(z))(T_z\Phi_g(v) \\ &\quad + \xi^1_{\mathcal{J}^{-1}(\rho)}(z'), T_z\Phi_g(w) + \xi^2_{\mathcal{J}^{-1}(\rho)}(z')) \\ &= (\Phi_g^*\omega)(z)(v, w) + \omega(\Phi_g(z))(T_z\Phi_g(v), \xi^2_{\mathcal{J}^{-1}(\rho)}(z'))) \\ &\quad + \omega(\Phi_g(z))(\xi^1_{\mathcal{J}^{-1}(\rho)}(z'), T_z\Phi_g(w)) \\ &\quad + \omega(\Phi_g(z))(\xi^1_{\mathcal{J}^{-1}(\rho)}(z'), \xi^2_{\mathcal{J}^{-1}(\rho)}(z')).\end{aligned}$$

Since the $G$–action is canonical, we have $(\Phi_g^*\omega)(z)(v, w) = \omega(z)(v, w)$. Also, since $\xi^1_{\mathcal{J}^{-1}(\rho)}(z') \in E(z')$, there exists a $G$–invariant function $f \in C^\infty(M)^G$ such that $X_f(z') = \xi^1_{\mathcal{J}^{-1}(\rho)}(z')$. This allows us to write

$$\omega(\Phi_g(z))(\xi^1_{\mathcal{J}^{-1}(\rho)}(z'), \xi^2_{\mathcal{J}^{-1}(\rho)}(z')) = \mathbf{d}f(z') \cdot \xi^2_{\mathcal{J}^{-1}(\rho)}(z') = 0,$$

by the $G$–invariance of the function $f$. If $h \in C^\infty(M)^G$ satisfies $X_h(z') = T_z\Phi_g(v)$, it follows that

$$\omega(\Phi_g(z))(T_z\Phi_g(v), \xi^2_{\mathcal{J}^{-1}(\rho)}(z'))) = \mathbf{d}h(z') \cdot \xi^2_{\mathcal{J}^{-1}(\rho)}(z')) = 0.$$

Analogously, we can conclude that $\omega(\Phi_g(z))(\xi^1_{\mathcal{J}^{-1}(\rho)}(z'), T_z\Phi_g(w)) = 0$ and hence

$$\omega_\rho([z']_\rho)([v']_\rho, [w']_\rho) = \omega(z)(v, w) = \omega_\rho([z]_\rho)([v]_\rho, [w]_\rho),$$

which guarantees that (4.4) is a good definition of $\omega_\rho$, consistent with (4.3). Notice that $\omega_\rho$ is smooth since $\pi_\rho^*\omega_\rho$ is smooth and it is also closed since $\omega$ is closed and $\pi_\rho$ is a surjective submersion.

We now show that $\omega_\rho$ is non degenerate, which concludes the proof. Let $[z]_\rho \in M_\rho$ and $[v]_\rho \in T_{[z]_\rho}M_\rho$ be such that for all $[w]_\rho \in T_{[z]_\rho}M_\rho$

$$\omega_\rho([z]_\rho)([v]_\rho, [w]_\rho) = 0.$$

which implies that $\omega(z)(v, w) = 0$ for all $w \in T_z(\mathcal{J}^{-1}(\rho)) = E(z)$. Let $f_1, f_2 \in C^\infty(M)^G$ be such that $v = X_{f_1}(z)$ and $w = X_{f_2}(z)$. Since $\omega(z)(v, w) = \mathbf{d}f_1(z) \cdot X_{f_2}(z) = 0$ for all $f_2 \in C^\infty(M)^G$, we conclude that $\mathbf{d}f_1(z) \in E(z)^\circ$. Let now $H := G_z$ and $L := N(H)/H$ (as usual, we will suppose for simplicity that $M_H$ is connected). Notice that the $L$–action on $M_H$ being canonical implies that for any $\eta \in \mathfrak{l}$ ($\mathfrak{l}$ denotes the Lie algebra of $L$), the vector field $\eta_{M_H}$ satisfies $\mathcal{L}_{\eta_{M_H}}\omega|_{M_H} = 0$ or, equivalently, $\eta_{M_H}$ is locally Hamiltonian. Let $\{\eta^1, \ldots, \eta^l\}$ be a basis of $\mathfrak{l}$ and $\mathbf{J}_L^{\eta^i}$ be the local Hamiltonian for $\eta^i_{M_H}$ defined in a neighborhood $U_i$ of $m$ in $M_H$. So, if $\eta = c_1\eta^1 + \cdots + c_l\eta^l$ is an arbitrary element of $\mathfrak{l}$, then $\eta_{M_H}$ is a locally Hamiltonian vector field with Hamiltonian function $\mathbf{J}_L^\eta = c_1\mathbf{J}_L^{\eta^1} + \cdots + c_l\mathbf{J}_L^{\eta^l}$, locally defined in $U := U_1 \cap \ldots \cap U_l$. We define $\mathbf{J}_L : U \to \mathfrak{l}^*$ by

$$\langle \mathbf{J}_L(z), \eta \rangle = \mathbf{J}_L^\eta(z).$$

It is easy to show that the map $\mathbf{J}_L$ has the following properties:

**(i)** $\ker(T_z\mathbf{J}_L) = (T_z(L \cdot z))^{\omega|M_H}$, for any $z \in U$.

**(ii)** $\operatorname{range}(T_z\mathbf{J}_L) = \mathfrak{l}^*$, for any $z \in U$.

**(iii)** Noether Theorem: The Hamiltonian flow associated to $L$–invariant functions in $M_H$ leaves the connected components of the level sets of $\mathbf{J}_L$ invariant.

Now, using the properties of $\mathbf{J}_L$ and the fact that $\mathbf{d}f_1(z) \in E(z)^\circ$, we can write

$$v = X_{f_1}(z)$$
$$= B^\sharp(z)(\mathbf{d}f_1(z)) \in B^\sharp(z)(E(z)^\circ) \cap E(z) \subset B^\sharp(z)((\ker(T_z\mathbf{J}_L))^\circ) \cap T_z M_H.$$

At the same time,

$$B^\sharp(z)((\ker(T_z\mathbf{J}_L))^\circ) \cap T_z M_H$$
$$= (\ker(T_z\mathbf{J}_L))^\omega \cap T_z M_H = T_z(L \cdot z) \subset T_z(G \cdot z)$$

Therefore $v = X_{f_1}(z) \in T_z(G \cdot z) \cap E(z) = T_z(G_\rho \cdot z) = \ker(T_z\pi_\rho)$, by Proposition 4.5. Consequently $[v]_\rho = 0$, as required. ∎

**4.8 Theorem** (Optimal Reduction of Hamiltonian dynamics).  *Consider a symplectic manifold $(M, \omega)$. Let $G$ be a Lie group acting properly and canonically on $M$. Let $E$ be the $G$–characteristic distribution with optimal momentum map $\mathcal{J} : M \to M/G_E$ and $h \in C^\infty(M)^G$ be a $G$–invariant Hamiltonian. Then*

(i) *The Hamiltonian flow $F_t$ of $X_h$ leaves the level sets $\mathcal{J}^{-1}(\rho)$ of $\mathcal{J}$ invariant and commutes with the $G$–action, hence if $\rho \in M/G_E$ satisfies the closedness hypothesis, $F_t$ induces a flow $F_t^\rho$ on $M_\rho$, uniquely determined by*

$$\pi_\rho \circ F_t \circ i_\rho = F_t^\rho \circ \pi_\rho, \tag{4.5}$$

*where $i_\rho : \mathcal{J}^{-1}(\rho) \hookrightarrow M$ is the canonical injection and $\pi_\rho : \mathcal{J}^{-1}(\rho) \to M_\rho$ is the projection.*

(ii) *The flow $F_t^\rho$ is Hamiltonian in $(M_\rho, \omega_\rho)$, with Hamiltonian function $h_\rho \in C^\infty(M_\rho)$ defined by*

$$h_\rho \circ \pi_\rho = h \circ i_\rho.$$

*We will call $h_\rho$ the **reduced Hamiltonian**. The vector fields $X_h$ and $X_{h_\rho}$ are $\pi_\rho$–related.*

(iii) *Let $k \in C^\infty(M)^G$ be another $G$–invariant function. Then, $\{h, k\}$ is also $G$–invariant and $\{h, k\}_\rho = \{h_\rho, k_\rho\}_{M_\rho}$, where $\{\cdot, \cdot\}_{M_\rho}$ denotes the Poisson bracket on $M$ associated to the symplectic structure $\omega_\rho$.*

**Proof.** Part (i) is a consequence of the Optimal Noether Theorem and the $G$–equivariance of the flow $F_t$. Parts (ii) and (iii) are a straightforward verification that uses the $G$–invariance of $h$ and the definition of $\omega_\rho$ given by expression (4.3). ∎

## 4.3 Comparison of the optimal and the Marsden–Weinstein reductions

Suppose now that the $G$–action in the statement of Theorem 4.7 is globally Hamiltonian, that is, there exists a globally defined equivariant momentum map $\mathbf{J}: M \to \mathfrak{g}^*$ that allows us to perform symplectic reduction in the spirit of Marsden and Weinstein [1974]. What is the relation between the Marsden–Weinstein reduced spaces obtained via the use of $\mathbf{J}$ and the optimal reduced spaces constructed using Theorem 4.7?

In Theorem 3.6 we showed that in the globally Hamiltonian case $\mathcal{J}^{-1}(\rho) = \mathbf{J}^{-1}(\mu) \cap M_H$, where $\mu = \mathbf{J}(m)$ and $H = G_m$, for some $m \in \mathcal{J}^{-1}(\rho)$ (all along this section we will assume that $\mathbf{J}^{-1}(\mu) \cap M_H$ is connected). Also, it is easy to show that in that situation $G_\rho = N_{G_\mu}(H) = N_G(H) \cap G_\mu$, with $G_\mu$ the coadjoint isotropy of $\mu \in \mathfrak{g}^*$. Indeed if $g \in G_\rho$ then $g \cdot \rho = \rho$. By the definition of the $G$–action on $M/G_E$, this implies that $\mathcal{J}(g \cdot m) = \mathcal{J}(m)$, or equivalently, both $g \cdot m$ and $m$ are in the same $G_E$–orbit, that is, there is a $G$–equivariant element $\mathcal{F}_T \in G_E$ such that $\mathcal{F}_T(m) = g \cdot m$. The $G$–equivariance of $\mathcal{F}_T$ implies that $m$ and $g \cdot m$ have the same isotropy subgroup, hence

$$H = G_m = G_{g \cdot m} = gG_m g^{-1} = gHg^{-1},$$

which implies that $g \in N(H)$. At the same time, Noether's Theorem for $\mathbf{J}$, as well as its $G$–equivariance implies that

$$g \cdot \mu = g \cdot \mathbf{J}(m) = \mathbf{J}(g \cdot m) = \mathbf{J}(\mathcal{F}_T(m)) = \mathbf{J}(m) = \mu,$$

which implies that $g \in G_\mu$. We therefore have that $G_\rho \subset N_{G_\mu}(H)$. The reverse inclusion is trivial once we assume the connectedness of $\mathbf{J}^{-1}(\mu) \cap M_H$. In conclusion, we have that

$$M_\rho = \mathbf{J}^{-1}(\mu) \cap M_H/N_{G_\mu}(H) = \mathbf{J}^{-1}(\mu) \cap M_H/(N_{G_\mu}(H)/H). \qquad (4.6)$$

When there are no singularities, the isotropies of all the elements in $\mathcal{J}^{-1}(\rho)$ are trivial ($H = \{e\}$) and therefore $M_H = M$ and $N_{G_\mu}(H) = G_\mu$. Consequently, in this case $M_\rho = \mathbf{J}^{-1}(\mu)/G_\mu$ and **the optimal and Marsden–Weinstein reduced spaces coincide**.

In the presence of singularities, the optimal reduced spaces (4.6) **coincide with the singular reduced spaces** introduced in Sjamaar and Lerman [1991] (for compact groups at zero momentum) and in Ortega [1998];

Ortega and Ratiu [2002] (for proper actions at arbitrary momentum values). Indeed,

$$M_\rho = \mathbf{J}^{-1}(\mu) \cap M_H/(N_{G_\mu}(H)/H) \simeq (\mathbf{J}^{-1}(\mu) \cap M_{(H)}^{G_\mu})/G_\mu,$$

where $M_{(H)}^{G_\mu} = \{m \in M \mid G_m = gHg^{-1}, g \in G_\mu\}$. See Sjamaar and Lerman [1991]; Ortega [1998]; Ortega and Ratiu [2002] for the details.

**Acknowledgments:** We are grateful for useful conversations with A. Alekseev, A. Blaom, V. Ginzburg, A. Knutson, B. Kostant, R. Loja Fernandes, J. Marsden, and A. Weinstein. This research was partially supported by the European Commission and the Swiss Federal Government through funding for the Research Training Network *Mechanics and Symmetry in Europe* (MASIE). Tudor Ratiu's Research was also partially supported by the Swiss National Science Foundation.

34     J-P Ortega and T. S. Ratiu